\documentstyle[txmac,a4,amssymb,
twoside,%
nocaphead,%
epsf,%
varthm,%
rotate,%
myrot,%
mypic,times,mathptm]{article}

\advance\oddsidemargin by -0.7cm
\advance\evensidemargin by -0.7cm
\advance\textwidth by 1.4cm

\def\mynewtheo#1#2{%
\newtheorem{@#1}{#2}[section]%
\newenvironment{#1}{\begin{@#1}\rm}{\end{@#1}}}

\mynewtheo{lem}{Lemma}
\mynewtheo{exer}{Exercise}
\mynewtheo{theo}{Theorem}
\mynewtheo{rem}{Remark}
\mynewtheo{defi}{Definition}
\mynewtheo{conj}{Conjecture}
\mynewtheo{corr}{Corollary}
\mynewtheo{prop}{Proposition}
\mynewtheo{question}{Question}
\mynewtheo{exam}{Example}

\newenvironment{eqn}{\begin{equation}}{\end{equation}}

\newenvironment{theorem}{\begin{theo}}{\end{theo}}

\begin{document}

\makeatletter

\newenvironment{myeqn*}[1]{\begingroup\def\@eqnnum{\reset@font\rm#1}%
\xdef\@tempk{\arabic{equation}}\begin{equation}\edef\@currentlabel{#1}}
{\end{equation}\endgroup\setcounter{equation}{\@tempk}\ignorespaces}

\newenvironment{myeqn}[1]{\begingroup\let\eq@num\@eqnnum
\def\@eqnnum{\bgroup\let\r@fn\normalcolor 
\def\normalcolor####1(####2){\r@fn####1#1}%
\eq@num\egroup}%
\xdef\@tempk{\arabic{equation}}\begin{equation}\edef\@currentlabel{#1}}
{\end{equation}\endgroup\setcounter{equation}{\@tempk}\ignorespaces}

\newcount\case@cnt
\newenvironment{caselist}{\case@cnt0\relax}{}

\newcommand{\mybin}[2]{\text{$\Bigl(\begin{array}{@{}c@{}}#1\\#2%
\end{array}\Bigr)$}}
\def\overtwo#1{\mbox{\small$\mybin{#1}{2}$}}
\newcommand{\mybr}[2]{\text{$\Bigl\lfloor\mbox{%
\small$\displaystyle\frac{#1}{#2}$}\Bigr\rfloor$}}
\def\mybrtwo#1{\mbox{\mybr{#1}{2}}}

\def\case{\advance\case@cnt by 1\relax{\bf\relax
\the\case@cnt\relax. Case.\ \ignorespaces}%
\edef\@currentlabel{\the\case@cnt}}

\def\myfrac#1#2{\raisebox{0.2em}{\small$#1$}\!/\!\raisebox{-0.2em}{\small$#2$}}
\def\sbx#1{\mbox{\scriptsize$\ds#1$}}

\def\mmyfrac#1#2{\raisebox{0.2em}{\scriptsize$#1$}\!/\!\raisebox
{-0.2em}{\scriptsize$#2$}}

\def\myeqnlabel{\bgroup\@ifnextchar[{\@maketheeq}{\immediate
\stepcounter{equation}\@myeqnlabel}}

\def\@maketheeq[#1]{\def\theequation{#1}\@myeqnlabel}

\def\@myeqnlabel#1{%
{\edef\@currentlabel{\theequation}
\label{#1}\enspace\eqref{#1}}\egroup}

\def\rottab#1#2#3{
\expandafter\advance\csname c@#3\endcsname by -1\relax
\centerline{%
\rbox{\centerline{\vbox{\setbox1=\hbox{#1}%
\hbox to \wd1{\hfill\vbox{{%
\caption{#2}}}\hfill}%
\vskip6mm
\box1}}%
}
}%
}

\def\epsfs#1#2{{\epsfxsize#1\relax\epsffile{#2.eps}}}

{\def\thefootnote{\fnsymbol{footnote}}
\footnotetext[1]{Supported by a DFG postdoc grant.} 

\author{A. Stoimenow\footnotemark[1]\\[2mm]
\small Department of Mathematics, \\
\small University of Toronto,\\
\small Canada M5S 3G3\\
\small e-mail: {\tt stoimeno@math.toronto.edu}\\
\small WWW: {\hbox{\tt http://www.math.toronto.edu/stoimeno/}}
}}

\title{\large\bf \uppercase{%
Generating functions, Fibonacci numbers}\\[2mm]
\uppercase{and rational knots}\\[4mm]
{\small\it This is a preprint. I would be grateful
for any comments and corrections!}%
}

\date{\normalsize Current version: \curv\ \ \ First version:
\makedate{13}{11}{1999}}


\maketitle

\makeatletter

\def\nowns{\not\owns}
\def\ex{\exists\,}
\def\fa{\forall\,}
\let\nb\nabla
\let\reference\ref
\let\ay\asymp
\let\pa\partial
\let\ap\alpha
\let\be\beta
\let\zt\zeta
\let\Gm\Gamma
\let\gm\gamma
\let\de\delta
\let\dl\delta
\let\Dl\Delta
\let\eps\epsilon
\let\lm\lambda
\let\Lm\Lambda
\let\sg\sigma
\let\vp\varphi
\let\om\omega

\let\es\enspace
\let\wh\widehat
\let\sm\setminus
\let\tl\tilde
\def\tr{\tl r}
\def\tq{\tl q}
\def\dt{\det}
\def\sgn{\mathop {\operator@font sgn}}
\def\rk{\mathop {\operator@font rank}}
\def\ncap{\not\mathrel{\cap}}
\def\cf{\text{\rm cf}\,}
\def\md{\max\deg}
\def\mc{\max\cf}
\def\Ra{\Rightarrow}
\def\lra{\longrightarrow}
\def\lmt{\longmapsto}
\def\mt{\mapsto}
\def\fra{\leftrightarrow}
\def\ul{\underline}
\def\so{\Rightarrow}
\def\So{\Longrightarrow}
\let\ds\displaystyle
\let\ol\overline
\def\li{\int\limits}

\def\ffrac#1#2{\mbox{\small$\ds\frac{#1}{#2}$}}

\newbox\@tempboxb
\def\sfrac#1#2{\setbox\@tempboxa=\hbox{$#2$}
\@tempdima=\ht\@tempboxa 
\advance\@tempdima by 2pt\relax 
\setbox\@tempboxb=\hbox{\vbox to \@tempdima{\vss}\unhbox\@tempboxa}
\frac{#1}{\box\@tempboxb}}

\long\def\@makecaption#1#2{%
   \vskip 10pt
   {\let\label\@gobble
   \let\ignorespaces\@empty
   \xdef\@tempt{#2}%
   }%
   \ea\@ifempty\ea{\@tempt}{%
   \setbox\@tempboxa\hbox{%
      \fignr#1#2}%
      }{%
   \setbox\@tempboxa\hbox{%
      {\fignr#1:}\capt\ #2}%
      }%
   \ifdim \wd\@tempboxa >\captionwidth {%
      \rightskip=\@captionmargin\leftskip=\@captionmargin
      \unhbox\@tempboxa\par}%
   \else
      \hbox to\captionwidth{\hfil\box\@tempboxa\hfil}%
   \fi}%
\def\fignr{\small\sffamily\bfseries}%
\def\capt{\small\sffamily}%

\newdimen\@captionmargin\@captionmargin1cm\relax
\newdimen\captionwidth\captionwidth\hsize\relax

\def\eqref#1{(\protect\ref{#1})}

\def\proof{\@ifnextchar[{\@proof}{\@proof[\unskip]}}
\def\@proof[#1]{\noindent{\bf Proof #1.}\enspace}

\def\@mt#1{\ifmmode#1\else$#1$\fi}
\def\qed{\hfill\@mt{\Box}}
\def\qqed{\hfill\@mt{\Box\enspace\Box}}

\let\Bbb\bf

\def\cN{{\cal N}}
\def\cS{{\cal S}}
\def\cC{{\cal C}}
\def\cM{{\cal M}}
\def\cP{{\cal P}}
\def\tg{{\tilde g}}
\def\tZ{{\tilde Z}}
\def\fg{{\frak g}}
\def\cZ{{\cal Z}}
\def\cD{{\cal D}}
\def\bR{{\Bbb R}}
\def\bC{{\Bbb C}}
\def\cE{{\cal E}}
\def\bZ{{\Bbb Z}}
\def\bN{{\Bbb N}}
\def\bQ{{\Bbb Q}}
\def\QI{{\Bbb Q}\cup\{\infty\}}
\def\RI{{\Bbb R}\cup\{\infty\}}

\def\bysame{\same[\kern2cm]\,}

\def\br#1{\left\lfloor#1\right\rfloor}
\def\BR#1{\left\lceil#1\right\rceil}

\def\TM{$^\text{\raisebox{-0.2em}{${}^\text{TM}$}}$}

\def\abstractname{}

{\let\@noitemerr\relax
\vskip-2.7em\kern0pt\begin{abstract}
\noindent{\bf Abstract.}\enspace
We describe rational knots with any of the possible
combinations of the properties (a)chirality, (non-)positivity,
(non-)fiberedness, and unknotting number one (or higher),
and determine exactly their number for a given number of
crossings in terms of their generating functions. 
We show in particular how Fibonacci numbers occur in the
enumeration of fibered achiral and unknotting number one
rational knots. Then we show how to enumerate rational
knots by crossing number and genus and/or signature. This
allows to determine the distribution of these invariants
among rational knots. We give also an application to
the enumeration of lens spaces.\\[1mm]
\end{abstract}
}

{\parskip0.2mm\tableofcontents}
\vspace{7mm}

\section{\label{sect1}Introduction}

A natural question one can ask in knot theory is how many
different knots or links, possibly of some special class,
there are of given crossing number
(that is, minimal number of crossings in any of their
diagrams). Clearly, to have a satisfactory approach to
such a problem, a good understanding of the class in
question is necessary. For an arbitrary knot or link,
the problems to identify it from a given diagram, and (hence also)
to determine its crossing number, although solved in theory
by Haken \cite{Haken}, are impracticably complicated. Thus an even
approximate enumeration of general knots and links by
crossing number seems so far impossible. However, some
bounds are known. In \cite{Welsh}, Welsh proved that
this number is exponentially bounded in the crossing
number $n$, with an upper bound to the base of the exponential
of $13.5$.

Even if a class of links is well-understood, still
its enumeration may be difficult. An example of such a class
are the prime alternating links. Such links have been classified
(contrarily to Haken, in a very practicable manner) in
\cite{MenThis}, and the determination of their crossing number
was settled in \cite{Kauffman,Murasugi,This} (both results
having been conjectured decades before by Tait).
Thus one can algorithmically generate the table of
links of certain (not too high) crossing number $n$, and hence
in particular determine (by ``brute force'') how many of them
there are \cite{HTW}. However, a reasonable expression for the
numbers thus obtained is not known, and possibly does not exist.
Only recently, Sundberg and Thistlethwaite \cite{SunThis} obtained
asymptotical estimates, accurate up to a linear factor in $n$.
(This slight inexactness was later removed in a note of Schaeffer
and Kunz-Jacques \cite{SC}.) In particular, they determined the
base of the exponential growth of these numbers to about $6.14$.
There has been other recent work \cite{ZZ}, which exhibits
a deep connection to statistical mechanics. This approach,
however, even if more effective than brute force enumeration,
is still very involved, and not yet made mathematically rigorous.
In \cite{nlpol}, I used the Sundberg-Thistlethwaite method
to improve Welsh's upper bound on the rate of growth of the number of
arbitrary links to about $10.3$.

Using quite different methods, basing on the
theory of Wicks forms, in joint work with A.~Vdovina
\cite{SV}, we determined the asymptotical behaviour of
the number of alternating knots of given genus up to a scalar
(depending on the genus).

All of these results are asymptotical and do not give exact formulas.
The only so far known such formulas concern the special class
of rational knots. In \cite{ErnSum}, Ernst and Sumners gave formulas
for the exact number of arbitrary and achiral rational knots and links
of given crossing number. The method they applied is again different
from the previously mentioned, and bases on Schubert's classification
\cite{Schubert} in terms of iterated (or continued) fractions.

In this paper, we will refine the results of Ernst
and Sumners for knots by considering three further
properties: positivity, fiberedness, and unknotting number one.
Together with achirality, these four properties subdivide the
class of rational knots into 16 subclasses, given by demanding or
excluding any of the properties. Only some of these subclasses are easy
to understand, since the properties defining them are causally dependent
(for example, positivity and achirality are mutually exclusive). Still
many of the classes are non-trivial, and apparently nothing about
them was so far known. We will obtain a description of all of
these classes, which allows to find an exact formula
for their size by crossing number. (In case there are only
few knots in the class, we will give them directly.) Usually,
it will be most convenient to give the numbers by means of their
generating functions, which turn out to be all rational functions
(a new, rather unexpected, justification for the designation of these
knots as rational).

Most interestingly, two of our enumeration problems
turn out to be directly related to Fibonacci numbers.
This way we have the possibly first explicit appearance of
this common integer sequence in a knot theoretically related
enumeration problem. A previous good candidate for such a problem
was the dimension of the space of primitive Vassiliev knot invariants
by degree. The apparent occurrence of Fibonacci numbers therein 
originally led to some excitement, until computer calculation
\cite{BarNatanVI} gave a disappointing result in degree 8,
where the dimension in question was 12, and not 13. (Now this
problem is known to be extremely hard and, if at all, will unlikely
offer such an elegant solution, see \cite{ChmDuz,Zagier}.)

We start with some preliminaries in \S\reference{Sp}, which
occupies the rest of this section, containing standard
definitions, facts, and conventions. The enumeration results
concerning knots with the aforementioned four properties will
be discussed then in \S\reference{Su1}--\reference{sP}.
Our method will be to study the effect of (combinations of)
these properties on the form of the iterated
fractions associated to the rational knots. It will be
in particular decisive to understand, how the two normal
forms, of all integers positive, and of all integers even,
transform into each other. While the description of fibered and
achiral rational knots in terms of their iterated fraction
is classical, the property of unknotting number one has
been made very approachable only by the more recent work
of Kanenobu and Murakami \cite{KanMur}. For positive rational
knots a convenient description will have to be worked out below.

One of the two enumeration results involving Fibonacci numbers will be
presented here first only as an inequality. We will remark that
the other (reverse) inequality depends on the truth of a
conjecture of Bleiler for fibered rational knots (henceforth
considered and meant only for unknotting number one). Namely,
in \cite{Bleiler} he conjectured that any rational knot
realizes its unknotting number in a rational diagram corresponding to
the expression of its iterated fraction with all integers even. This
conjecture was disproved by Kanenobu and Murakami \cite{KanMur},
quoting the counterexample $8_{14}$, which is not fibered. Since
Bleiler's conjecture now again turns out to be relevant in the
fibered case, it will be the matter of new consideration.

In a note \cite{rat2}, written after this paper was begun, but
already published, we announce and complete the results of
this paper, by giving a first proof of the ``fibered'' Bleiler
conjecture. There we also formulate a statement about
unimodular matrices which is related to this conjecture.
For the historical reason to explain a part of the result of
\cite{rat2}, we include here a section \S\reference{SBl}.
In this section we establish a relation between
the conjecture a products of certain unimodular matrices,
and use this relation to obtain some results related to it.
Before this, in \S\reference{Sff}, we discuss how to enumerate
counterexamples to Bleiler's conjecture (of which the
Kanenobu-Murakami knot $8_{14}$ is the simplest one),
by classifying their even-integer notations. This
also leads to a new proof of (a generalization
of) the ``fibered'' Bleiler conjecture. 

In \S\reference{S6.1} and \S\reference{Ssig}, we give a few other
formulas, 
including one determining the number of rational knots
of given genus or signature. The formulas arising here contain
several variables and are much more involved. Some of them are not
rational, but all can be given in closed form. They yield by
substitutions the Ernst--Sumners result, and also several formulas
obtained previously in this paper. We apply an integration
method allowing to build the generating function of the product
of two sequences and to ``select'' certain parts of a
multivariable generating function.

The final enumeration results will concern lens spaces
by fundamental group, by using their correspondence to
rational knots of given determinant, of which they are
the 2-fold branched coverings. In the enumeration some
exceptional (duplication) series of determinants occur,
and the question whether they intersect non-trivially
is related to the integer solutions of a certain
hyperelliptic equation. 

We conclude the introduction with a remark addressing rational links.
We decided to leave them completely out of the discussion in this
paper. One reason is that there will be already enough to say
on knots. Secondly, at least most of the arguments can
be adapted to links. (In fact, links occur naturally jointly
with knots at some places, and we will have then to
artificially get disposed of them.) However, for links
also unpleasant questions connected with
orientation come in, and would make the approach more
technical than methodical.

\subsection{Preliminaries and notations\label{Sp}}

The Fibonacci numbers $F_n$ are a very popular integer
sequence. These numbers can be defined recursively by
$F_0=F_1=1$ and $F_n=F_{n-1}+F_{n-2}$, explicitly by
\[
F_{n-1}=\frac{1}{\sqrt{5}}\left[\left(\frac{1+\sqrt{5}}{2}\right)^n-
\left(\frac{1-\sqrt{5}}{2}\right)^n\right]\,,
\]
and also by the generating function
\[
\sum_{n=0}^\infty F_nx^n=\frac{1}{1-x-x^2}\,.
\]
See your favorite calculus textbook, or
\cite[sequence A000045]{Sloane} for an extensive compilation
of references. Due to this simple property Fibonacci numbers
appear very often in many unrelated situations and it is
always amazing to see them come up in some mathematical problem.

$\#S=|S|$ are alternative designations for the cardinality of $S$.

A \em{knot} is a $C^1$ embedding $K:\,S^1\hookrightarrow S^3$ (for
convenience henceforth identified with its image) up to isotopy.
Usually knots are represented by diagrams,
plane curves (images of $K$ under the
projection of $\bR^3=S^3\sm\{*\}$ onto a generic hyperplane) with
transverse self-intersections (crossings) and distinguished (over)%
crossing strand (a connected component of the preimage under $K$ of
a neighborhood of the crossing).

A knot $K$ is called \em{fibered}, if $S^3\sm K$ is a bundle over $S^1$
with fiber being a Seifert surface for $K$, an embedded in $S^3$
punctured compact orientable surface $S$ with $\pa S=K$
(see \cite{Gabai}).
In this case, by the theorem of Neuwirth-Stallings, $S$ has minimal
genus among all Seifert surfaces for $K$ (called the \em{genus} 
$g(K)$ of $K$) and is unique up to isotopy.

A knot has \em{unknotting number} one if it
has some diagram, such that a crossing change $
\diag{0.7em}{2}{2}{\piclinewidth{20}\picmultiline{-8.5 1 -1 0}{0 0}{2 2}
\picmultiline{-8.5 1 -1 0}{2 0}{0 2}}\,\to\,
\diag{0.7em}{2}{2}{\piclinewidth{20}\picmultiline{-8.5 1 -1 0}{2 0}{0 2}
\picmultiline{-8.5 1 -1 0}{0 0}{2 2}}\,$, creates (a diagram of)
the unknot (the knot with diagram $\diag{0.7em}{2}{2}{
\piclinewidth{20}\piccircle{1 1}{1}{}}$). More generally,
one defines the unknotting number $u(K)$ of a knot $K$ as the minimal
number of crossing changes in any diagram of $K$ needed to turn
$K$ into the unknot (see e.g. \cite{KanMur,Lickorish}).

The \em{crossing number} of a knot
is the minimal crossing number of all its diagrams.

A knot $K$ is called \em{achiral} (or amphicheiral) if there exists an
isotopy turning it into its mirror image in $S^3$, otherwise $K$ is
called \em{chiral}.

The writhe is a number ($\pm1$), assigned to any crossing in a link
diagram. A crossing as in figure \ref{figwr}(a) has writhe 1 and
is called positive. A crossing as in figure \ref{figwr}(b) has writhe 
$-1$ and is called negative.

\begin{figure}[htb]
\[
\begin{array}{c@{\qquad}c}
\diag{6mm}{1}{1}{
\picmultivecline{0.18 1 -1.0 0}{1 0}{0 1}
\picmultivecline{0.18 1 -1.0 0}{0 0}{1 1}
} &
\diag{6mm}{1}{1}{
\picmultivecline{0.18 1 -1 0}{0 0}{1 1}
\picmultivecline{0.18 1 -1 0}{1 0}{0 1}
}
\\[2mm]
(a) & (b)
\end{array}
\]
\caption{\label{figwr}}
\end{figure}

A knot is called \em{positive} if it has a positive diagram,
i.~e. a diagram with all crossings positive. See for example
\cite{Cromwell,MorCro,positive}.

The \em{braid index} of a knot is the minimal number of
strands of a braid which closes up to the knot; see
\cite{Murasugi2}. In \cite{Murasugi3}, one finds a
definition and properties of the \em{signature}.

A knot $K$ is \em{rational} (or $2$-bridge), if it has bridge
number 2, where the bridge number is half of the smallest number of
critical points of a Morse function on $K$. In \cite{Schubert},
rational knots have been classified by the iterated fractions
corresponding to their \em{Conway notation} \cite{Conway}.
See Goldman and Kauffman \cite{GK} for a more modern account.

Let the \em{iterated fraction} $[[s_1,\dots,s_m]]$ for
integers $s_i$ be defined inductively by $[[s]]=s$ and
\[
[[s_1,s_2,\dots]]=s_1+\frac{1}{[[s_2,\dots]]}\,.
\]
Note: there is another convention of building iterated fractions,
in which the `$+$' above is replaced by a `$-$'. See e.g.
\cite{Murasugi2}. Latter is more natural in some sense
(see the proof of theorem \reference{thls}), but the
permanent sign switch makes it (at least for me) more
unpleasant to work with in practice. Thus we stick to
the version with `$+$'.

The rational knot or link $S(p,q)$ in Schubert's \cite{Schubert}
notation has the Conway \cite{Conway} notation $c_{n}\ c_{n-1}\ \dots
\ c_{1}$, when the $c_i$ are chosen so that 
\begin{eqn}\label{ci}
[[c_{1},c_{2},c_{3},\dots,c_n]]=\frac{p}{q}\,.
\end{eqn}
Since when replacing integers with variables the
Conway notation $a_1\,a_2\,\dots a_n$ in its original form
becomes somewhat illegible, we will sometimes put the sequence into
parentheses, or also use the alternative
designation $C(a_1,\dots,a_n)$ for this sequence.
Thus $S(p,q)=C(c_{n},c_{n-1},\dots,c_{1})$.
We also abbreviate repeating subsequences as powers,
for example $(4(12)^21^33)=(412121113)$. We call
the numbers $c_i$ also \em{Conway coefficients} of the notation.

Note that $S(-p,-q)$ is the same knot or link as $S(p,q)$, while
$S(-p,q)=S(p,-q)$ is its mirror image. $S(p,q)$ is a knot for
$p$ odd and a 2-component link for $p$ even. The number $p$
is the \em{determinant} of $K$, given by $|\Dl_K(-1)|$,
where $\Dl$ is the Alexander polynomial \cite{Alexander}.

Without loss of generality one can assume that $(p,q)=1$,
$|q|<|p|$, and that (exactly) one of $p$ and $q$ is even.
(If both $p$ and $q$ are odd, we replace $q$ by
$q\pm |p|$, the sign being determined by the
condition $|q|<|p|$.)

Then we can choose all $c_i$ in \eqref{ci} to be even (and non-zero).
It is known that, with this choice of $c_i$, their number
$n=2g(S(p,q))$ is equal to twice the genus of $S(p,q)$ for $p$ odd
(i.e. a rational knot). To fix a possible ambiguity between a
diagram and its mirror image, we consider the crossings corresponding
to the entry $c_i$ in the Conway notation to have writhe
$(-1)^{i-1}\sgn(c_i)$.

For the purpose of calculating with iterated fractions,
it will be helpful to extend the operations `$+$' and `$1/.$' to
$\QI$ by $1/0=\infty,\ 1/\infty=0,\ k+\infty=\infty$ for any $k\in
\bQ$. The reader may think of $\infty$ as the fraction $1/0$, to which
one applies the usual rules of fraction arithmetics and reducing.
In particular reducing tells that $-1/0=1/0$ so that for us
$-\infty=\infty$. This may appear at first glance strange, but has 
a natural interpretation in the rational tangle context.

It will be helpful to introduce some notation for subsequences
of the sequence of integers giving the Conway notation for some
rational knot. We most commonly denote such subsequences by letters
towards the end of the alphabet like $x$ or $y$, while single
integers will be called $a,b,\dots$. Define for a finite
sequence of integers $x=(a_1,\dots,a_n)$ its \em{reversion} (or \em{%
transposition}) $\ol x:=(a_n,\dots,a_1)$ and its \em{negation}
by $-x:=(-a_1,\dots,-a_n)$. If $x=\pm \ol x$, we call $x$
\em{(anti)palindromic}. For $y=(b_1,\dots,b_m)$
the term $x,y$ denotes the concatenation of both
sequences $(a_1,\dots,a_n,b_1,\dots,b_m)$. Similarly one
defines concatenation with a single integer, for example 
$x,b=(a_1,\dots,a_n,b)$ etc.

\begin{figure}[htb]
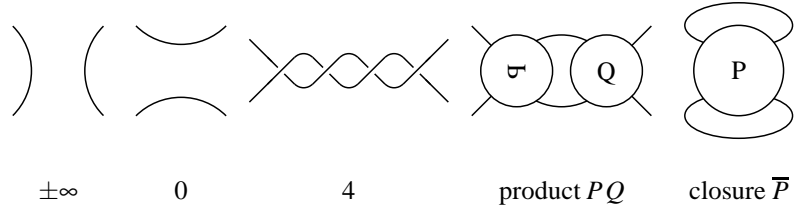

\[
\begin{array}{*5c}
\diag{6mm}{2}{2}{
  \piccirclearc{-1 1}{1.41}{-45 45}
  \piccirclearc{3 1}{1.41}{135 -135}
} & 
\diag{6mm}{2}{2}{
  \piccirclearc{1 -1}{1.41}{45 135}
  \piccirclearc{1 3}{1.41}{-135 -45}
} &
\diag{6mm}{4}{1}{
  \picmultigraphics{4}{1 0}{
    \picline{0.2 0.8}{0.8 0.2}
    \picmultiline{-6 1 -1.0 0}{0.2 0.2}{0.8 0.8}
  }
  \picmultigraphics{3}{1 0}{
    \piccirclearc{1 0.6}{0.28}{45 135}
    \piccirclearc{1 0.4}{0.28}{-135 -45}
  }
  \picline{-0.2 -0.2}{0.3 0.3}
  \picline{-0.2 1.2}{0.3 0.7}
  \picline{4.2 -0.2}{3.7 0.3}
  \picline{4.2 1.2}{3.7 0.7}
} &
\diag{6mm}{4}{2}{
  \picline{0 0}{0.5 0.5}
  \picline{0 2}{0.5 1.5}
  \picline{4 0}{3.5 0.5}
  \picline{4 2}{3.5 1.5}
  \piccirclearc{2 0.5}{1.3}{45 135}
  \piccirclearc{2 1.5}{1.3}{-135 -45}
  \pictranslate{1 1}{
    \picrotate{90}{
      \picscale{-1 1}{
	\picfilledcircle{0 0}{0.8}{P}
      }
    }
  }
  \picfilledcircle{3 1}{0.8}{Q}
} &
\diag{6mm}{2.4}{3}{
  \picmultigraphics{2}{0 2}{
    \picellipse{1.2 0.5}{1.2 0.5}{}
  }
  \picfilledcircle{1.2 1.5}{1}{P}
} \\[12mm]
\pm\infty & 0 & 4 & 
\text{ product $P\,Q$ } & \text{ closure $\ol{P}$ }
\end{array}
\]
\caption{Conway's primitive tangles and tangle operations.\label{figt}}
\end{figure}

Figure \ref{figt} shows how to obtain a diagram of the rational knot
or link from its Conway notation. It is the closure of
the rational tangle with the same notation. The convention in
composing the tangles is that a Conway notation with no
negative integers gives an alternating diagram. For more details see
\cite[\S 2.3]{Adams}.

Since by \cite{Kauffman,Murasugi,This} (reduced) alternating
diagrams have minimal crossing number, the crossing number
of a rational knot is the sum of the integers in its
Conway notation with all integers positive.

%
%
%

A good reference on generating function theory is \cite{Wilf}.

Before we start with our results, we make two remarks.

First, we will adopt the convention of considering rational
knots \em{up to mirroring}. This is, a rational knot and its
mirror image will be considered equivalent, and hence counted once.
Since we will have formulas for the number of rational knots with
specific properties (counted up to mirroring) and for
the number of achiral rational knots with the same specific
properties, one can easily obtain from both numbers the number of
rational knots with the same properties with chiral knots and their
mirror images counted separately. An exception to this convention will
be \S\reference{Ssig} and \S\reference{Sfg}, where the sensitivity of
the signature
under mirroring forces care to be taken. There we will specify
in each statement whether we count chiral pairs once or twice.

When considering rational knots up to mirroring, the Conway notation
with all integers positive is determined up to reversal and the
up to the ambiguity $\dots,n-1,1 \fra \dots,n$ at the end.
The Conway notation with all integers even and non-zero
is determined up to reversal and the simultaneous negation
of all entries. In both cases the reversal of notation
corresponds to the identity $S(p,q)=S(p,\pm q^{-1})$.
Here $q^{-1}$ is the multiplicative inverse of $q$ in
$\bZ_p^*$, the group of units of $\bZ_p=\bZ/p\bZ$, and the
sign is positive or negative depending on whether the
Conway notation has odd or even length.

Also, we will content ourselves only with interesting combinations of
the four properties. For some of the remaining combinations,
the results are known, sometimes even in greater generality than
just for rational knots. We refer to \cite{gen1,positive} for
the treatment of these cases, and do not consider them here. (One
could certainly prove some of these results also from our setting,
but such an attempt does not seem any longer relevant.) For
other combinations of properties, the description easily follows
from what we will prove below. For some of them, a subset of
the properties already restricts sufficiently the class, and
the check of the remaining properties on the few knots is
easy. The remaining enumerations follow by simple inclusion-exclusion
arguments. In such cases we mostly waive on presenting the
results explicitly here and leave them to the reader.

\section{Rational knots of unknotting number one\label{Su1}}

We start with the description and enumeration of rational knots of
unknotting number one for given crossing number. Here, unlike
in subsequent sections, we first use the notation with positive
Conway coefficients. (We will return to unknotting number one later,
when armed with a more effective method.)

\begin{theorem}
If $K$ is a rational knot of unknotting number one, then it has a
Conway notation with all coefficients positive, which is in at least
one of the types listed below. (In the first five cases the entry `$-1$'
indicates the crossing to be switched to unknot the knot.)
{\nopagebreak
\def\labelenumi{\roman{enumi})}\mbox{}\\[-18pt]
\def\theenumi{\roman{enumi})}
\begin{enumerate}
\item\label{type1} $x,n,-1,1,n-1,\bar x,c$
\item\label{type2} $x,n-1,1,-1,n,\bar x,c$
\item\label{type3} $a,x,n,-1,1,n-1,\bar x,a\pm 1$
\item\label{type4} $n+1,-1,1,n-1,1,c$ and $n-1,1,-1,n-1,1,c$
(degenerate cases of \ref{type1} and \ref{type2})
\item\label{type5} $n-1,1,-1,n\pm 1$ (degenerate case of \ref{type3})
\item\label{type6} $2,n$
\end{enumerate}
}
Here $x$ denotes a (possibly empty) sequence of positive integers,
and $a$, $c$ and $n$ are positive integers, so that all entries in
the above sequences (except the `$-1$') are positive. Also, unlike
elsewhere, $x$ is considered up to the ambiguity $n,\dots=1,n-1,\dots$
for $n>1$. (Thus for example
the sequence $(5,2,3,1,1,4,2,4,1,7)$ is considered of type
\ref{type2} with $c=7$, $n=4$ and $x=(5,2)=(1,4,2)$.)
\end{theorem}

\proof It was proved in \cite[\S 3.1]{rtan}, that a rational
knot of unknotting number one has an alternating diagram 
of unknotting number one, and hence all alternating diagrams
have this property. Consider the alternating diagram of the Conway
notation with all integers $a_1,\dots,a_k$ positive. If the (unknotting)
crossing change occurs in a group of $\ge 2$ half-twists,
then the only such case is \reference{type6}. Else we need
to switch `$1$'$\to$`$-1$'. In this case after this change we obtain
modulo mirroring a closed rational tangle with iterated fraction $1/n$
for some $n\in\bN$. Modulo transposition of the notation, we may
assume $n=\pm 1$ (case \reference{type3}) or that the (sub)tangle
with Conway notation $a_1,\dots,a_{k-1}$ turns into the
0-tangle under the crossing change (giving cases \reference{type1}
and \reference{type2} with $c=a_k$). The almost-symmetry in the
first three cases arises 
when analyzing the iterated fraction from left and right until
the crossing changed. Up to a correction $n,\dots\to1,n-1,\dots$
in their inner ends, and the ambiguity $\dots,p=\dots,p-1,1$
at the outer ends (because only their iterated fraction is relevant)
they must be transposed. This explains the occurrence of
$x$ and $\bar x$. (The ambiguity at the outer end of $\bar x$
changes the knot if not at outermost position in the notation.)
The degenerate cases \reference{type4} and \reference{type5} occur
when the fraction expression has length one. \qed

From the theorem (and the lack of essential restrictions to $x$)
the enumeration of unknotting number one rational knots of given
crossing number is straightforward (but rather tedious by virtue of
having to take care of duplicatedly counted cases and the
ambiguity for $x$). Thus it is clear how to obtain the following
corollary, which was suggested empirically. However, instead of
going now into unpleasant details, we will later give a much
more elegant proof.

\begin{corr}\label{ttt}
If $c_n$ denotes the number of rational unknotting
number one knots of $c$ crossings (chiral pairs counted only once),
then these numbers are given basically by powers of 2, namely
via the generating (rational) function
\[
\sum_{n=1}^\infty c_nx^n\,=\,x^3+x^4(x+1)\left[\frac{2}{1-
2x^2}+\frac{1}{x^2-1}\right]+\frac{x^8}{x^4-1}\,.
\]
In particular, $\lim\limits_{n\to\infty}\sqrt[n]{c_n}=\sqrt{2}$.
\end{corr}

It is worth mentioning that for every fourth crossing number the
number of rational unknotting number one knots does not increase
compared to the next crossing number~-- this is possibly not what
one may expect!

\begin{corr}\label{crau1}
The number of achiral unknotting number one
rational knots of $c$ crossings is $2$ for $c=10+6k$, $k\ge 0$, and $1$
for other even $c\ge 4$. More exactly, these knots are those with
Conway notation $(n11n)$ and $(3(12)^k1^4(21)^k3)$.
\end{corr}

\proof It is known that $C(a_1,\dots,a_n)$ with all $a_i>0$ is achiral
iff the sequence $a_1,\dots,a_n$ is (up to the ambiguity $\dots,n-1,1
\fra \dots,n$) palindromic of even length. The result then is a direct 
verification of the palindromicity of the patterns of the above $6$
cases. The series $(n11n)$ clearly comes from case \reference{type5},
while $(3(12)^k1^4(21)^k3)$ for $k>0$ is in case \reference{type1}
and for $k=0$ in case \reference{type4}. The other cases only give
at best alternative representations for $4_1=(22)$ and $6_3=(2112)$.
\qed

\begin{corr}
Except for the trefoil and figure eight knot,
all unknotting number one rational knots have in their alternating
diagrams at most two crossings, such that switching any single 
one of them unknots the knot.
\end{corr}

\proof The theorem shows that if a `$1$' is changed to `$-1$' to
unknot, then the number of integers left and right from it differs by
at most three. This leaves at most 4 (neighbored) positions. The
degenerate cases are easily excluded, and considering \reference{type1},
\reference{type2} and \reference{type3}, one finds that only the edge
`$1$' in a subsequence of `$1$'s can unknot, and at most one of these
edge `$1$'s does, if the sequence is of length two (except for case
\reference{type5}). \qed

Clearly the knots where (exactly) two such crossings exist
include the achiral ones given in corollary \reference{crau1}.
We leave it to the reader to modify the proof of corollary
\reference{crau1} and to show that the remaining knots are of the
forms $(32^k1^32^{k+1})$ and $(2^k1^32^{k})$. (This result
was again suggested by computer calculation, and I have not
carried out a rigorous proof.)

\section{Fibered rational knots}

For the following results it is more convenient to work with
the (unique up to reversal and negation) expression of the iterated
fraction by even (non-zero) integers rather than natural numbers.
(The number of all these even integers is always even and equal to
the double genus of the knot.) The key point is how to extract
the crossing number out of this representation. The result
is given in the following lemma, which will be of central
importance throughout the rest of the paper.

\begin{lem}\label{lerr}
If $a_1,\dots,a_{2g}$ are even (non-zero) integers, then the
crossing number of $C(a_1,\dots,a_{2g})$ is 
\[
\sum_{i=1}^{2g}|a_i|-\,\#\{\,1\le i<2g\,:\,a_ia_{i+1}<0\,\}\,.
\]
(In fact, the formula still holds if all $|a_i|\ge 2$ not necessarily
even.)
\end{lem}

\proof
We remarked that the crossing number result for alternating diagrams
\cite{Kauffman,Murasugi,This} implies that the crossing number
of a rational knot is the sum of the integers in its
Conway notation with all integers positive. Thus we need to
account for the change of the sum of the $|a_i|$, when
transforming the Conway notation with all integers even into
the one with all integers positive. This is a repeated application
of the iterated fraction identity $[[x,a,-b,y]]=[[x,a-1,1,b-1,-y]]$
(with $x$ and $y$ subsequences and $a$ and $b$ integers).
The claim then follows by induction on the number of
such applications needed. \qed

\begin{theorem}
If $c_n$ denotes the number of fibered rational
knots of $n$ crossings (chiral pairs counted only once),
then these numbers are given by the generating function
\[
\sum_{n=1}^\infty c_nx^n\,=\,-\frac{x^3(1+x)(x^4+x^3+x^2-1)}
{(x^4+2x^3+x^2-1)(x^4+x^2-1)}\,.
\]
In particular,
$\ds\lim\limits_{n\to\infty}\sqrt[n]{c_n}=\frac{1+\sqrt{5}}{2}$.
\end{theorem}

The proof is a prototype of argument that will occur in several
more complex variations later.

\proof A rational knot is fibered iff all even integers $a_i$ in
its iterated fraction expression are $\pm 2$. This is a well-known
fact which seems to have been (algebraically) noted explicitly
in this form first by Lines and Weber \cite{LW}, although it is also 
a consequence of the (much older) result of Murasugi \cite{Murasugi4},
as we shall briefly argue. A geometric proof can be also given,
for example using the method of \cite{Gabai}.

The diagram of the closure of a rational tangle with all integers
even is the Murasugi sum of connected sums of Hopf bands with $a_i/2$
full twists. Thus from \cite{Murasugi4} the multiplicativity of
the leading coefficient $\mc\Dl$ of the Alexander polynomial
under Murasugi sum implies
\begin{eqn}\label{pp}
\mc\Dl_{C(a_1,\dots,a_{2g})}\,=\,\pm 2^{-2g}\,\prod_{i=1}^{2g}\,a_i\,.
\end{eqn}
If the knot is fibered, $\mc\Dl=\pm 1$, and hence all $a_i=\pm 2$.
Contrarily, if all $a_i=\pm 2$, the knot has a surface which
is a plumbing of Hopf bands with one full twist each, and hence a
fiber surface.

In the case a rational knot is fibered, each $a_i\pm 2$, except
the first one, according to lemma \reference{lerr},
contributes one to the crossing number of the knot, if it
follows a $\pm 2$ of the different sign, and two otherwise.
Thus, by ignoring the contribution of the first $\pm 2$,
we are left by counting compositions into parts $1$ and $2$
of $n-2$ of odd length up to transposition. (A composition
of a certain number is writing it as a sum of numbers, whose
order is relevant.)

To pass from this to the generating function of the theorem is
a matter of some combinatorial calculation. One uses the generating
function
\[
f_1(x)\,=\,\frac{x+x^2}{1-\bigl(x+x^2\bigr)^2}\,,
\]
for the number of odd length compositions into parts $1$ and $2$,
and
\[
f_2(x)\,=\,\frac{x+x^2}{1-x^2-x^4}
\]
for the number of palindromic ones.

If we fix the first number $a_1=2$ up to mirroring, the notations
define the same knot iff they differ by transposition and possible
negation (so as the initial term to become positive).

In this situation, $f_1$ counts all knots we like by their notations
twice, except the ones with palindromic and antipalindromic notations.
These are enumerated exactly by $f_2$. To see this, one needs to
remark that the sequence of $1$'s and $2$'s contributing from each $a_i$
to the crossing number, with the initial $2$ coming from $a_1$ omitted,
is palindromic iff the Conway notation made up of the $a_i$ (without the
initial one $a_1$ omitted) is palindromic or antipalindromic.

Thus the generating function we seek is simply $(f_1+f_2)/2$. \qed

\begin{rem}
One can, of course, give using partial fraction decomposition
an explicit formula for the $c_n$ in terms of (negative powers of)
the zeros of the denominator polynomial of the generating function,
from which the limit property (that is, the justification to write above
`$\lim$' rather than `$\limsup$') follows, but the resulting expression
should be less pleasant, so we waive on its derivation.
\end{rem}

\begin{rem}
One could, in a similar way, show that the number of rational knots $K$
with $\mc\Dl_K$ being up to sign a fixed natural number $n$ give the
Taylor coefficients of a rational function. The complexity of this
function will roughly depend on the complexity of the prime
decomposition of $n$. This relies on the fact that for $a_1,\dots,a_{2g}
$ even (and non-zero), we have the relation \eqref{pp}.
\end{rem}

The fact that we count compositions into parts $1$ and $2$
already suggests the relation to Fibonacci numbers.
Now comes the enumeration result where they appear explicitly.

\begin{theorem}
The number of rational fibered achiral knots of $n$ crossings is
$F_{n/2-2}$ for $n$ even (and $0$ for $n$ odd).
\end{theorem}

\proof It is known, see \cite{Schubert}, that a rational knot is achiral
iff its Conway notation with all integers even is palindromic.
Then it follows directly from the lemma that the crossing number must
be even (this follows more generally for alternating knots from
\cite{Kauffman,Murasugi,This}). Considering only the first
half of the (palindromic) sequence, we see that again
we count compositions into parts $1$ and $2$, this time of $\myfrac{1}
{2}(n-4)$, but neither the restriction on the number of parts
(genus $-1$), nor the factoring out of transpositions are necessary.
Thus the result follows. \qed

The Fibonacci numbers also occur when considering
rational fibered unknotting number one knots.

\begin{theorem}\label{thFb}
The number of rational fibered unknotting number one knots of $n$
crossings is $2F_{\br{n/2-3}}$ for $n\ge 6$.
\end{theorem}

As most of the results
before, theorem \reference{thFb} was also suggested by computer,
which calculated the various sequences above up to 26 crossings,
and confirmed for this sequence equality with the doubled Fibonacci
numbers. The completion of the proof depends on the truth of the
``fibered'' Bleiler conjecture. This is motivation enough to come
back to this problem in more detail towards the end of the paper.

\proof[of theorem \reference{thFb}]
Consider the diagrams $C(\pm 2,x,\pm 2,-\bar x)$, $x$ being
a sequence of $\pm 2$'s and $-\bar x$ its negated transposed. 
Clearly such knots unknot by switching a crossing counted
by the $\pm2$ in between $x$ and $-\bar x$. It is also easy to see
that if a diagram $C(b_1,\dots,b_{2g})$ with all $b_i$ even
and non-zero is to be unknottable by one crossing change,
it must be of this form. The only possible cancellations
near a zero entry are of the form
\begin{eqn}\label{yyz}
(\dots,a,0,b,\dots)=(\dots,a+b,\dots)\,,
\end{eqn}
and when only non-zero entries remain,
the notation does not represent the unknot.

Thus from now on consider Conway notations of the form
\begin{eqn}\label{ppp}
C(a_0,a_1,\dots,a_k,\pm 2,-a_k,\dots,-a_1)\,,
\end{eqn}
with all $a_i=\pm 2$. Our concern will be to
count such notations by crossing number, as given in the
lemma \reference{lerr}.

Now, for a given knot, the Conway notation with all numbers even is
unique up to negating all numbers and transposition. In order
to avoid duplicate counting, we must take care what notations still
fit into the form \eqref{ppp} after some of these transformations.

Clearly, negating all numbers preserves the form \eqref{ppp},
but to get disposed of this transformation, we can simply declare that
we count only forms with $a_0>0$.

Then we must find out which sequences
\eqref{ppp} remain of this from after
transposition. For such sequences one sees that the first and last
entries determine the rest of the sequence.
Since we restricted ourselves only to sequences with $a_0>0$,
we see that demanding $a_0=-a_1=2$ forces the sequence to
become palindromic, and hence it is not counted twice.
(This sequence then corresponds to the knots with notation
$(n11n)$ given in corollary \reference{crau1}.)

Now we can apply the previous arguments. Again one counts
compositions into parts $1$ and $2$ coming from the subsequence
$x=(a_1,\dots,a_k)$, and the equality of the numbers one obtains
for $n$ and $n+1$ if $n$ is even comes from the switch of
signs in $x$ together with the sign of the middle $\pm 2$. 
Switching just the sign of the middle $\pm 2$, fixing $x$,
accounts for the factor $2$. 


We proved so far that there are \em{at least} as many knots
as we claimed in the formulation of theorem \reference{thFb}.
To remove that `at least', one needs that any fibered rational knot of
unknotting number one should realize its unknotting number in a rational
diagram of all Conway coefficients even. This was conjectured for
arbitrary rational knots by Bleiler \cite{Bleiler}, but disproved by
Kanenobu and Murakami \cite{KanMur}, quoting the counterexample
$8_{14}$ (which, however, is not fibered). Thus, the confirmation
of Bleiler's conjecture for fibered rational knots and unknotting
number one is equivalent to establishing equality in (and completing
the proof of) the above theorem.
%
As noted, the statement we require was proved in \cite{rat2}, but
another and more generalized proof (which will also lead to
generalizations of this theorem) will be given in \S\reference{Sff}.
\qed

\begin{rem}
To describe the fibered rational knots which are both of unknotting
number one and achiral, one uses corollary \reference{crau1}. The knots
in the first family there are fibered (they are closed alternating
3-braids), while those in the second family are not. To see latter fact,
the reader may convince himself, that the crossings counted by
the initial and terminal `$3$' in the Conway notation correspond
to reverse(ly oriented) half-twists:
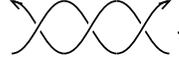
\[
\diag{7mm}{3}{2}{
  \picPSgraphics{0 setlinecap}
  \pictranslate{1 1}{
    \picrotate{-90}{
      \rbraid{0 -0.5}{1 1}
      \rbraid{0 0.5}{1 1}
      \rbraid{0 1.5}{1 1}
      \pictranslate{-0.5 0}{
      \picvecline{0.03 1.95}{0 2}
      \picvecline{0.03 -.95}{0 -1}
    }
  }
 }
}\,.
\]
It follows from the description of $\mc\Dl$
on alternating diagrams given in \cite{Cromwell} that an alternating
diagram with $\ge 3$ reverse half-twist crossings always has
$|\mc\Dl|>1$, and hence never represents a fibered link.
\end{rem}

\begin{table}[ptb]
\captionwidth\vsize\relax
\newpage
\vbox to \textheight{\vfil
\rottab{ 
\hbox to \vsize{\footnotesize\hss
\begin{mytab}{|c||*{24}{r|}}%
  { & \multicolumn{24}{|r|}{ } }
  \hline [2mm]%
$n$ & $3$ &   $4$ &   $5$ &   $6$ &   $7$ &   $8$ &    $9$ &    $10$ &   $11$ &   $12$ &    $13$ &    $14$ &    $15$ &     $16$ &  $17$ &     $18$ &      $19$ &      $20$ &      $21$ &      $22$ &       $23$ &       $24$ &   $25$ &        $26$ \\[2mm]%
\hline
\hline [2mm]%
f & $1$ & $1$ & $1$ & $2$ & $3$ & $4$ & $7$ & $10$ & $16$ & $25$ & $40$ & $62$ & $101$ & $159$ & $257$ & $410$ & $663$ & $1062$ & $1719$ & $2764$ & $4472$ & $7209$ & $11664$ & $18828$ \\
fa & & $1$ & & $1$ & & $2$ & & $3$ & & $5$ & & $8$ & & $13$ & & $21$ & & $34$ & & $55$ & & $89$ & & $144$ \\
u & $1$ & $1$ & $1$ & $3$ & $3$ & $6$ & $7$ & $15$ & $15$ & $30$ & $31$ & $63$ & $63$ & $126$ & $127$ & $255$ & $255$ & $510$ & $511$ & $1023$ & $1023$ & $2046$ & $2047$ & $4095$ \\
au & & $1$ & & $1$ & & $1$ & & $2$ & & $1$ & & $1$ & & $2$ & & $1$ & & $1$ & & $2$ & & $1$ & & $1$ \\
fu & $1$ & $1$ & & $2$ & $2$ & $2$ & $2$ & $4$ & $4$ & $6$ & $6$ & $10$ & $10$ & $16$ & $16$ & $26$ & $26$ & $42$ & $42$ & $68$ & $68$ & $110$ & $110$ & $178$ \\
p & $1$ & & $2$ & & $5$ & & $12$ & & $30$ & & $76$ & & $195$ & & $504$ & & $1309$ & & $3410$ & & $8900$ & & $23256$ & \\
$\sg0$ & & $1$ & & $3$ & $2$ & $9$ & $6$ & $29$ & $30$ & $99$ & $112$ & $351$ & $450$ & $1275$ & $1734$ & $4707$ & $6762$ & $17577$ & $26208$ & $66197$ & $101862$ & $250953$ & $395804$ & $956385$
\\[2mm]
\hline
\end{mytab}
\hss}
}{The number of rational knots of $n \le 26$ crossings with some
combinations of the properties \ul achirality,
\ul unknotting number one, \ul positivity, and \ul fiberedness
(the combination is indicated by joining the initials
of the properties considered). The last line contains
the number of rational knots of zero signature, whose
determination will be explained later.\protect\label{tab1}}{table}
\vss}
\newpage
\end{table}

\section{Positive rational knots\label{sP}}

Positive knots (see \cite{Cromwell,MorCro,positive,Yokota} for example)
have been around for a while in knot theory, but apparently no
special attention was given to the rational ones among them.
We start by a description of such knots, again
using the expression with all Conway coefficients even.

\begin{lem}\label{rpos}
If $a_1,\dots,a_{2g}$ are even integers, then the rational knot
$C(a_1,\dots,a_{2g})$ is positive, iff all $a_i$ alternate in sign, i.e.
$\#\{\,1\le i<2g\,:\,a_ia_{i+1}<0\,\}=2g-1$\,.
\end{lem}

\proof If a rational (or alternating) knot is positive, then by
\cite{restr,Nakamura} so is any of its alternating diagrams, and
hence by \cite{Murasugi3}, $\sg=2g$ (where $\sg$ is the \em{signature}
and $g$ the genus). If all $a_i$ alternate in sign, then
$C(a_1,\dots,a_{2g})$ is a positive diagram. However, is
some $a_i$ has the wrong sign, then $C(a_1,\dots,a_{2g})$ can
be obtained from a positive diagram by undoing positive/creating
negative reverse twists.
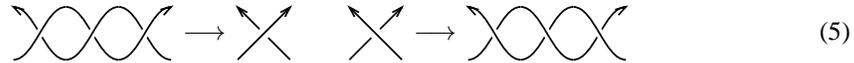
\begin{eqn}
\diag{7mm}{3}{2}{
  \picPSgraphics{0 setlinecap}
  \pictranslate{1 1}{
    \picrotate{-90}{
      \rbraid{0 -0.5}{1 1}
      \rbraid{0 0.5}{1 1}
      \rbraid{0 1.5}{1 1}
      \pictranslate{-0.5 0}{
      \picvecline{0.03 1.95}{0 2}
      \picvecline{0.03 -.95}{0 -1}
    }
    }
    }
} \,\lra\, 
\diag{7mm}{1}{1}{
    \picmultivecline{-7 1 -1.0 0}{1 0}{0 1}
    \picmultivecline{-7 1 -1.0 0}{0 0}{1 1}
} \qquad
\diag{7mm}{1}{1}{
    \picmultivecline{-7 1 -1.0 0}{0 0}{1 1}
    \picmultivecline{-7 1 -1.0 0}{1 0}{0 1}
} \,\lra\,
\diag{7mm}{3}{2}{
  \picPSgraphics{0 setlinecap}
  \pictranslate{1 1}{
    \picrotate{-90}{
      \lbraid{0 -0.5}{1 1}
      \lbraid{0 0.5}{1 1}
      \lbraid{0 1.5}{1 1}
      \pictranslate{-0.5 0}{
      \picvecline{0.03 1.95}{0 2}
      \picvecline{0.03 -.95}{0 -1}
    }
    }
    }
} 
\end{eqn}
Any of these moves does not augment $\sg$. Moreover, as in
this process at least once some $a_i=0$, giving a knot of smaller
genus, and as $\sg\le 2g$, $\sg$ strictly decreases. Then
$C(a_1,\dots,a_{2g})$ has $\sg<2g$, and the knot is not
positive. \qed

\begin{theorem}\label{thpo}
If $c_n$ denotes the number of rational positive
knots of $n$ crossings, then these numbers are given by the
generating function
\[
\sum_{n=1}^\infty c_nx^n\,=\,\frac{x^3-2x^5}{(1-3x^2+x^4)(1-x^2-x^4)}\,.
\]
In particular, all positive rational knots have odd crossing number, and
$\ds\lim\limits_{n\to\infty}\sqrt[2n+1]{c_{2n+1}}=\frac{1+\sqrt{5}}{2}$.
\end{theorem}

\proof Since now, by the lemma, if $C(a_1,\dots,a_{2g})$ is positive,
all $a_i$, $i>1$ contribute $a_i-1$ to the crossing number of $K$,
by artificially decreasing $a_1$ by $1$, we are left with counting
compositions of $n-1$ into (an even number of) odd parts up to
transposition.

Without factoring out transpositions, the number
of compositions of $n$ of this type is given by the generating
function
\[
\frac{1}{\sbx{1-\Bigl(\frac{\strut x}{1-x^2}\Bigr)^2}}\,.
\]
The number of palindromic compositions is the same as the
number of compositions of $n/2$ into (a now not necessarily even
number of) odd parts, whose generating function is
\[
\sfrac{1}{\sbx{1-\frac{x^2}{1-x^4}}}\,.
\]
Thus, accounting for the unknot, we have
\[
\sum_{n=1}^\infty c_nx^n\,=\,-x+\frac{x}{2}\left[
\frac{1}{\sbx{1-\Bigl(\frac{\strut x}{1-x^2}\Bigr)^2}}+
\sfrac{1}{\sbx{1-\frac{x^2}{1-x^4}}}\right]\,,
\]
whence the result. \qed

The theorem roughly suggests that there should be approximately
qualitatively equally many positive and fibered rational knots
up to a given crossing number. This should be contrasted to the
distribution of their iterated fractions: while $\{\,p/q\,:\,
S(p,q)\mbox{ is positive }\}$ should be dense in $\bR\sm(-1,1)$,
the closure of $\{\,p/q\,:\,S(p,q)\mbox{ is fibered }\}$ will have
dense \em{complement} (possibly even zero Lebesgue measure\footnote{%
According to a remark of A.~Sikora, not any complement of an
open dense subset must have zero measure. In fact, there are open
dense subsets in $\bR$ of arbitrarily small positive measure!}).

Of course, that positive knots have odd crossing number is not true
even for prime alternating knots, as shows $8_{15}$.

\begin{corr}\label{cff}
The only fibered positive rational knots are the $(2,n)$-torus
knots for $n$ odd.
\end{corr}

\proof Combining the lemma with the fiberedness property, we obtain
a notation $C((2,-2)^g)$, which belongs to the $(2,2g+1)$ torus knot.
\qed

This shows that almost all fibered positive knots are not rational.
In \cite{Murasugi2}, Murasugi mainly settled the problem of
non-alternation (so in particular non-rationality) for closed positive
braids. However, he needs the technical assumption that such knots
have positive braid representation of minimal strand number, and
moreover, not all  fibered positive knots are closures of positive
braids, as shows the example $10_{161}$ discussed in
\cite[example 4.2]{positive}. In \cite{restr}, we generalize
corollary \ref{cff} to alternating knots and links.

Table \reference{tab1} summarizes some of the numbers discussed above.


The previous arguments can be applied to several similar
enumeration problems. We discuss in some detail how to
obtain the number of rational knots of given genus and/or
given signature.

\section{Genus\label{S6.1}}

For the genus, one can prove

\begin{theorem}\label{thG}
If $c_{n,g}$ is the number of rational knots of $n$ crossings and genus $g$,
then
\begin{eqn}\label{cng}
f(x,z)\,=\,\sum_{g=1}^{\infty}\sum_{n=3}^{\infty}\,c_{n,g}x^nz^g\,=\,
-{\frac{{x^3}\,z\,\left( -1 + {x^3}\,z + {x^4}\,z + 
{x^2}\,\left( 1 + z \right)  \right) }{\left( 1 + x \right) \,
\left( 1 + {x^2} \right) \,
\left( -1 + 2\,x + {x^2}\,\left( -1 + z \right)  \right) \,
\left( -1 + {x^2}\,\left( 1 + z \right)  \right) }}
\end{eqn}
is a rational function in $x$ and $z$.
\end{theorem}

This is a similar, but slightly stronger, analogue of a result of \cite{gen1},
where we showed that
\[
f_g(x)\,=\,\sum_{n=3}^{\infty}\,c'_{n,g}x^n
\]
is a rational function in $x$ for fixed $g$, with $c'_{n,g}$ being the
number of \em{alternating} knots of $n$ crossings and genus $g$.
However, the dependence of $f_g$ on $g$ is too complicated to
let expect any nice (in particular, rational) expression for
the two-variable function \eqref{cng} in the alternating case.

\proof[of theorem \reference{thG}]
To prove the assertion for the genus, use that it is one half of
the number of entries in the even Conway notation. Let
$w=(a_1,\dots,a_{2g})$ be the sequence of these entries. Assume
$a_1>0$ to factor out one of the ambiguities. Define
as before a transformation of $w$ to a sequence $\hat w$ of positive
integers by $\hat w=(b_1,b_2,\dots,b_{2g})$, such that
\[
b_i\,=\,\left\{\,\begin{array}{cl}
|a_i| & \mbox{if $i=1$ or $a_{i-1}a_i>0$} \\
|a_i|-1 & \mbox{otherwise}
\end{array}\right.\,.
\]
Every sequence of positive integers with the first one even has
a unique preimage under $\wh{\,\cdot\,}$. We are interested in
counting those sequences $w$ such that the sum of entries of
$\hat w$ is $n$. Since $\wh{\,\cdot\,}$ is injective, this is
the same as counting compositions of $n$ into $2g$ (positive
integer) parts, the first one being even.
If $a_{n,g}$ is the number of such compositions,
then
\[
g(x,z)\,=\,\sum_{g=1}^{\infty}\sum_{n=1}^{\infty}\,a_{n,g}x^nz^g
\,=\,
\frac{x}{1+x}\left(\,\frac{1}{1-\ds\frac{zx^2}{(1-x)^2}}-1\,\right)\,.
\]

We count now every sequence $w$ once, but still there
are different sequences giving the same knot, coming
from the ambiguity of reversing the notation. Namely,
this always happens except if the sequence $w$ is
palindromic ($w=\bar w$) or anti-palindromic ($w=-\bar w$).
Let $y$ be the first half of $w$ (of length $g$). Then
$\hat y$ has a sum of entries either $n/2$ if $w$ is
palindromic, or $(n+1)/2$, if $w$ is anti-palindromic.

Thus for given $n$, only palindromic or only
anti-palindromic sequences $w$ occur, and their number
is the same as the number of compositions of $\BR{\ffrac{n}{2}}$
of length $g$ with the first integer being even.

If $b_{n,g}$ is the number of compositions
of $n$ of length $g$ with the first integer being even,
then
\[
h(x,z)\,=\,\sum_{g=1}^{\infty}\sum_{n=1}^{\infty}\,b_{n,g}x^nz^g
\,=\,
\frac{x}{1+x}\left(\,\frac{1}{1-\ds\frac{\strut zx}{1-x}}-1\,\right)\,.
\]
To replace $n$ by $\BR{\ffrac{n}{2}}$, one has to divide by $x$,
replace $x$ by $x^2$, and multiply by $x+x^2$.
\[
h_1(x,z)\,=\,
\frac{x+x^2}{1+x^2}\left(\,\frac{1}{1-\ds\frac{zx^2}{1-x^2}}-1\,\right)\,.
\]

Then $f(x,z)=\myfrac{1}{2}(g(x,z)+h_1(x,z))$. \qed


\begin{rem}\label{rem4}
Since the even-degree-$x$ part of $h_1$ counts the
palindromic sequences $w$, which correspond exactly to
the achiral knots, $\myfrac{1}{2}(h_1(x,z)+h_1(-x,z))$
gives the function enumerating achiral knots by genus.
\end{rem}

Using the result of Murasugi \cite[Theorem B (2)]{Murasugi2},
one can obtain a similar formula for counting by crossing number and
braid index:

\begin{theorem}
If $c_{n,b}$ is the number of rational knots of
$n$ crossings and braid index $b$, then
\[
\sum_{b=2}^{\infty}\sum_{n=3}^{\infty}\,c_{n,b}x^nz^b\,=\,
-\, \frac{x^3\,z^2\,\left( -1 - x\,z + 2\,x^4\,z^2 + x^5\,z^3 + 
x^2\,\left( 1 + z \right)  + x^3\,z\,\left( 2 + z \right)  \right) }
{\left( 1 + x \right) \,\left( -1 + x + 2\,x^2\,z \right) \,
\left( -1 + x^2 + 2\,x^4\,z^2 \right) }\,.
\]
\end{theorem}

\proof The proof is analogous and largely omitted. (Take, however,
care of the different convention for building iterated fractions.)
We remark only that instead of the previous function $g$ we must take
\[
 \frac{x\,z}{1 + x\,z}\,\left( \frac{1}
 {\ds 1 - \left[ {\left( 1 + \frac{1}{x\,z} \right) }\,
 {\left( \ds  \frac{1}{1 - x^2\,z} -1 \right) } \right]^2} -1 \right)\,,
\]
and instead of $h_1$
\begin{myeqn}{\qed}
\frac{x\,z\,\left( 1 + x\,z \right)}{1 + x^2\,z^2} \,
 \left( \sfrac{1} {\ds 1 - \left( 1 + \frac{1}{x^2\,z^2} \right) \,
 \left( \frac{1}{1 - x^4\,z^2} -1 \right) }^{\phantom{2} } -1 \right)\,.
\end{myeqn}

One can also count by genus and braid index \em{without} incorporating
the crossing number, as one observes that for given genus and
braid index there are only finitely many rational knots (a fact
which can be proved in larger generality). The discussion so far
should explain sufficiently how to proceed, so that
we leave this task to an interested reader.

\section{\label{Ssig}Signature}

Another variation of the enumeration problem (for which an analogue
for alternating knots, if it exists, is even harder to prove)
is to count rational knots by signature $\sg$.

One has the following formula for the signature (see
\cite[p.\ 71]{HNK}):

\begin{lem}\label{lmsg}
\begin{eqn}\label{opq}
\sg(C(a_1,\dots,a_{2g}))\,=\,\sum_{i=1}^{2g}\,(-1)^{i-1}\sgn(a_i)\,,
\end{eqn}
for $a_i\ne 0$ all even.
\end{lem}

This formula shows a close relationship between signature and
genus. Thus in this case we must again take care of the
genus, and so this is a refinement of the enumeration
by genus. Set in the sequel for simplicity $\chi'=1-\chi=2g$.

Since now $\sg$ depends on mirroring and takes negative values, we
must be careful about what and how exactly we are going to count.

There are several options how to avoid the chirality
and the negative value problems.

{\nopagebreak\parskip4pt 
\def\labelenumi{\arabic{enumi})}
\def\theenumi{\arabic{enumi})}
\mbox{}\\[-30pt]
\begin{enumerate}
\item\label{_type1} It appears suggestive to count
chiral pairs twice, since both knots give distinct contributions,
and there is no natural way to distinguish one of them. Then
we must deal with negative values of $\sg$, since
it is desirable to avoid negative powers in the generating series.
There are also two options:

{\nopagebreak
\def\theenumi{\arabic{enumi}}
\def\labelenumii{\alph{enumii})}
\def\theenumii{\alph{enumii})}
\mbox{}\\[-30pt]
\begin{enumerate}
\item\label{_typea} we count knots by $|\sg|\in[0,\chi']$, or
\item\label{_typeb} we count knots by $\sg+\chi'\in[0,2\chi']$.
\end{enumerate}
}

\item\label{_type2} Alternatively, but less naturally, one can declare
to count in a chiral pair only the knot with $\sg>0$, and if
both knots have $\sg=0$, to take any one of both, since their
contribution is the same. This has the advantage of also
eliminating the $\sg<0$ problem.
\end{enumerate}
}

We will thus for the rest of \S\reference{Ssig} specify according
to what convention we count rational knots by signature, and
in particular whether we count chiral pairs once or twice.

\subsection{$\sg$ with mirroring}

First we will deal with the version \reference{_typeb}.
It can be approached most naturally and leads to the ``simplest''
generating series. It fits into the picture we described
throughout the preceding discussion.

\begin{theorem}
Let $G_1$ be the function in 3 variables
$x$, $y$ and $z$ which counts in its Taylor coefficient of $x^my^lz^k$
the number of rational knots of crossing number $m$ with
$1-\chi=l$ and $1-\chi+\sg=k$, such that (unlike so far
in the paper) \em{both} knots in a chiral pair are counted. Then
$G_1$ is a certain rational function (shown in full form in figure
\reference{tb2}).
\end{theorem}

\begin{figure}[ptb]
\captionwidth0.98\vsize\relax
\newpage
\vbox to \vsize{\vfil
\rottab{%
\fbox{\parbox[t]{0.98\textheight}{%
\begin{eqnarray*}
G_1(x,y,z) & = &
  - x^3\,y^2\,\left( -1 - z^4 + 
         x^8\,z^2\,{\left( -1 + y^2\,z^2 \right) }^2\,
          \left( 1 + y^2\,z^2 \right)  - x\,\left( 1 + z^2 + z^4 \right)  + 
         x^6\,z^2\,\left( 1 + 2\,y^6\,z^6 + 2\,y^2\,\left( 1 + z^4 \right)  - 
            y^4\,z^2\,\left( 2 + 3\,z^2 + 2\,z^4 \right)  \right)  + 
	    \right. \\& &
         + x^7\,z^2\,\left( 1 + y^6\,z^6 + y^2\,\left( 1 + z^2 + z^4 \right)  - 
            y^4\,\left( z^2 + 3\,z^4 + z^6 \right)  \right)  + 
         x^2\,\left( -z^2 + y^2\,\left( 1 + 4\,z^4 + z^8 \right)  \right)  + 
         x^3\,\left( -z^2 + y^2\,
             \left( 1 + z^2 + 3\,z^4 + z^6 + z^8 \right)  \right)  + 
	     \\ & &
	     \left.  + x^5\,\left( 1 + z^2 + z^4 - 
            y^4\,z^4\,\left( 2 + z^2 + 2\,z^4 \right)  + 
            y^2\,\left( 1 + 2\,z^2 + 2\,z^6 + z^8 \right)  \right)  + 
         x^4\,\left( 1 - z^2 + z^4 - 3\,y^4\,\left( z^4 + z^8 \right)  + 
            y^2\,\left( 1 + 2\,z^2 + z^4 + 2\,z^6 + z^8 \right)  \right) 
         \right) \times \\ & & 
	 \times \frac{1}{\left( 1 + x \right) \,\left( 1 + x^2 \right) \,
       \left( 1 - x\,y\,\left( 1 + z^2 \right)  + 
         x^2\,\left( -1 + y^2\,z^2 \right)  \right) \,
       \left( 1 + x\,y\,\left( 1 + z^2 \right)  + 
         x^2\,\left( -1 + y^2\,z^2 \right)  \right) \,
       \left( 1 - x^2\,y^2\,\left( 1 + z^4 \right)  + 
         x^4\,\left( -1 + y^4\,z^4 \right)  \right) } \\[3mm]
\end{eqnarray*}
\begin{eqnarray*}
f_0(x) & = & 
     \frac{-x}{2\,\left( 1 + x \right) \,\left( 1 + x^2 \right) \,
     {\sqrt{\left( -1 + 4\,x^2 \right) \,\left( -1 + 4\,x^4 \right) }}}
  \,\left( -{\sqrt{1 - 4\,x^2}} + 
         2\,x^5\,{\sqrt{1 - 4\,x^2}} - {\sqrt{1 - 4\,x^4}} + 
         2\,{\sqrt{\left( -1 + 4\,x^2 \right) \,
              \left( -1 + 4\,x^4 \right) }}\, + 
	 \right. \\ & & +
         2\,x^4\,\left( {\sqrt{1 - 4\,x^2}} + {\sqrt{1 - 4\,x^4}} \right)  - 
         x\,\left( {\sqrt{1 - 4\,x^2}} + {\sqrt{1 - 4\,x^4}} - 
            2\,{\sqrt{\left( -1 + 4\,x^2 \right) \,
                 \left( -1 + 4\,x^4 \right) }} \right)  - 
		 \\ & & \left.
         - x^3\,\left( {\sqrt{1 - 4\,x^2}} + {\sqrt{1 - 4\,x^4}} - 
            2\,{\sqrt{\left( -1 + 4\,x^2 \right) \,
                 \left( -1 + 4\,x^4 \right) }} \right)  + 
         x^2\,\left( -{\sqrt{1 - 4\,x^2}} + {\sqrt{1 - 4\,x^4}} + 
            2\,{\sqrt{\left( -1 + 4\,x^2 \right) \,
                 \left( -1 + 4\,x^4 \right) }} \right)  \right)
\end{eqnarray*}
}}
}%
{\label{tb2}}{figure}
\vfil}
\newpage
\end{figure}

After our proof we will indicate how to proceed with enumeration
version \reference{_type2}, whose function can be expressed from the one
of version \reference{_typeb} by means of a certain complex integral
(and thus is no longer rational). The function for
enumeration version \reference{_typea} is obtained similarly,
and so we do not present it here.

\proof Let us start as before. Consider again a sequence $w$ of even
integers $w=(a_1,\dots,a_{2g})$ with $a_1>0$, and the associated
sequence $\hat w$. The formula for $\sg$ in the lemma can be
read as follows in terms of $\hat w$: 
subdivide $\hat w$ into subsequences starting with
an even integer, followed by some (possibly empty)
sequence of odd integers. Each such subsequence contributes
its length with alternating sign to the signature.
Call a subsequence $\sg$-positive or $\sg$-negative
dependingly on the sign of its contribution to $\sg$.

Let
%
%
\[
\hat F(x,y,z)=yz\frac{x^2}{1-x^2}\left(\frac{1}{1-\ds y\frac{
\strut xz}{1-x^2}}\right)\,.
\]
By the previous arguments we see that
\[
F_1(x,y)=\hat F(x,y,1)
\]
counts a single $\sg$-negative group
of entries by $\chi'$ in (powers of) $y$ and crossing number in $x$
(here $\chi'+\sg=0$). Similarly
\[
F_2(x,y,z)=\hat F(x,y,z^2)
\]
counts a single $\sg$-positive
group of entries by $\chi'$ in $y$, crossing number in $x$
and $\chi'+\sg$ in $z$.

Now $\hat w$ is made up of an arbitrary number of interchangingly
positive and negative subsequences, starting with a positive one.
Thus to count $\hat w$ we consider
\[
\r F(x,y,z)=\left(1+\frac{1}{F_2}\right)\frac{F_1F_2}{1-F_1F_2}\,,
\]
which counts an arbitrary sequence of $\sg$-positive/negative
groups by $\chi'$ in $y$ and crossing number in $x$.
This function now contains odd powers of $y$ ($=$values of $\chi'$).
They are discarded by setting
\[
F(x,y,z)=\frac{\r F(x,y,z)+\r F(x,-y,z)}{2}\,,
\]
which selects all knots ($1-\chi$ even), and
counts knots without factoring by palindromic ambiguity.

As before any knot, whose $w$ is not palindromic or antipalindromic,
is counted twice. However, here ``counted twice'' might have meant
that actually the knot and its mirror image have been counted,
thus contributing to two different coefficients in the power series.

Keeping in mind this point, we turn to 
care about palindromic sequences.
{\nopagebreak
\def\labelenumi{\arabic{enumi})}\mbox{}\\[-18pt]
\def\theenumi{\arabic{enumi})}
\begin{enumerate}
\item Consider the antipalindromic case.
$w$ is automatically of even length. Let $w'$ be the first
half of $w$. To simplify notation, let
\[
c(w)\,=\,|\hat w|_1=\sum b_i,\quad\chi'(w)=\mbox{length of }w,\quad
\sg(w)\,=\,\sum(-1)^{i-1}\sgn(a_i)\,.
\]
We remarked in the genus enumeration that
\[
c(w')=\frac{c(w)+1}{2} \quad\mbox{and}\quad
\chi'(w')=\frac{\chi'(w)}{2}\,.
\]
It remains to observe that also 
\[
\sg(w')=\frac{\sg(w)}{2}\,,
\]
which easily follows from the definition.

Thus antipalindromic cases are counted by
\[
F_0(x,y,z)=\frac{1}{x}\r F(x^2,y^2,z^2)\,.
\]

\item In the palindromic case we may have $\chi'$ odd.
However, it is easy to see that $\chi'(w)$ is odd if and only if
$c(w)$ is so, so that working only with even powers of $x$
will ensure that we count only knots. Assuming $\chi'(w)$ is
even and letting $w'$ be the first
half of $w$, we have
\[
c(w')=\frac{c(w)}{2},\quad \chi(w')=\frac{\chi'(w)}{2}
\quad\mbox{and}\quad \sg(w)=0\,,
\]
so that
\[
(\chi'+\sg)(w)=\chi'(w)=2\chi'(w')\,.
\]
Then we obtain $F_3$ enumerating palindromic cases from $\r F$
by replacing $y$ with $y^2z^2$ and $z$ by $1$, as $\sg(w')$
has no contribution to $\sg(w)$:
%
\[
F_3(x,y,z)=\r F(x^2,y^2z^2,1)\,.
\]
%
\end{enumerate}
}

Let
\[
G(x,y,z)=\frac{F(x,y,z)+F_0(x,y,z)+F_3(x,y,z)}{2}\,.
\]
Now the coefficient of $x^ky^{2g}z^l$ + the coefficient of
$x^ky^{2g}z^{2g-l}$ in $G(x,y,z)$ counts the number of rational knots
with crossing number $k$, genus $g$ and $2g+\sg=l$, where for each
chiral pair either only one knot is recorded, or both
are recorded with factor $\myfrac{1}{2}$
(the coefficients of $G$ lie only in $\bZ\cup\bZ+\myfrac{1}{2}$ !).
$F_3$ counts the achiral ones.

To count for each chiral pair both knots,
%
%
%
we set
\[
G_1(x,y,z)=G(x,y,z)+G(x,yz^2,1/z)-F_3(x,y,z)\,,
\]
which counts both knots in chiral pairs by $\chi'$ and $\chi'+\sg$
(the variable substitution in the second term accounts for
$y^lz^{l+k}\lra y^lz^{l-k}$). Thus $G_1$ is the function we sought.
\qed

\begin{rem}
One has the ($\sg$-forgetting) identity
\[
G_1(x,y,1)\,=\,2f(x,y^2)-\frac{h_1(-x,y^2)+h_1(x,y^2)}{2}\,,
\]
with $f$ being the 2-variable function in theorem \reference{thG},
and $h_1$ the one occurring in its proof. See remark \reference{rem4}.
Also, it is easy to see from the proof of lemma \reference{rpos},
that $G_1(x,1,0)$ enumerates negative
rational knots by crossing number. Since they correspond bijectively
to positive knots, $G_1(x,1,0)$ must coincide with the
function we obtained in theorem \reference{thpo}. Both
identities are easily verified.
\end{rem}

\subsection{$|\sg|$ without mirroring}

In $G_1$ a knot and its mirror image are represented by two
coefficients, for $\pm|\sg|$, i.e. for monomials
$y^{2g}z^{2g\pm|\sg|}$, which accounts for the
symmetry of $G_1$ under $(y,z)\to (yz^2,1/z)$. We can eliminate
this redundancy and count rational knots according to version
\reference{_type2}. Then we have

\begin{prop}\label{prJ}
Let $J$ be the function in $x$, $y$ and $z$ which counts in its
coefficient of $x^my^lz^k$ rational knots of crossing number $m$ with
$1-\chi=l$ and $|\sg|=k$, such that again only one knot 
in a chiral pair is counted. Then $J$ is a certain closedly
expressible function (too complicated to display).
\end{prop}

\proof To obtain $J$ from $G_1$,
basically we want to ``cut off'' terms in $G_1$ of monomials $y^lz^k$
with $k<l=2g$ (so far $[G_1]_{y^lz^kx^m}\ne 0$ for $0\le k\le 2l$),
and substitute $y^lz^k\lra y^lz^{k-l}$. We must care
about the chiral knots with $\sg=0$. Thus we consider 
\[
G_2(x,y,z)=G_1(x,y,z)+F_3(x,y,z)\,,
\]
and must multiply the coefficients in $G_2$ of $x^my^lz^k$ by
\[
\left\{\begin{array}{cc}0 & k<l\\\myfrac{1}{2}& k=l\\1& k>l
\end{array}\right.\,,
\]
and make the variable substitution $y\to y/z$.

If $H=\sum a_ix^i$ and $G=\sum b_ix^i$ converge in
a complex neighborhood of $0$, then for any $\ap\in [0,1]$ and
$|x|$ small
\begin{eqn}\label{intf}
\sum a_ib_ix^i\,=\,\li_0^1\,G(x^\ap e^{2\pi it})
H(x^{1-\ap}e^{-2\pi it})\,dt\,.
\end{eqn}
This formula is justified under the assumption of absolute convergence
and integrability of the limit function. The values of $x$, for which
this happens usually depends on $\ap$, but it is only
important that it contains a set with
a convergence point. Then, if the integral can be solved
in closed form for these $x$, by the uniqueness of the holomorphic
extension it also holds for all $x$ for which the series on
the left converges.

With this formula (under the convergence and integrability assumption,
which can be achieved with $\ap=0$ for $|y|,|z|<1$ and $|x|<1/2$), the
function $J(x,y,z)$ we seek can be expressed by an integral
\begin{eqn}\label{inl}
J(x,y,z)\,=\,\frac{1}{2\pi}\li_0^{2\pi}
G_2(x,ye^{-is},e^{is})\left(\frac{1}{1-ze^{-is}}-\frac{1}{2}\right)\,ds
\,.
\end{eqn}
This integral is, regrettably, too hard to solve pleasantly even with
the help of a computer, not from the point of view of method, but
of the structural complexity of the expressions to handle.
As we stated the proposition only qualitatively, we mostly avoid the
quotation of exact calculation results.

The integral was evaluated as follows. First, one uses
standard substitution $t=\tan \myfrac{s}{2}$, with which it
turns into a rational integral
\[
\,\li_{-\infty}^{\infty}\,G_2\left(x,y\frac{1-t^2-2it}{1+t^2},
\frac{1-t^2+2it}{1+t^2}\right)\,\left(\frac{1+t^2}{1+t^2-z(
1-t^2-2it)}-\frac{1}{2}\right)\,\frac{dt}{\pi\,(1+t^2)}\,.
\]
This integral can be solved by calculating the residues of the
(meromorphic) integrand in the upper half-plane. One integrates
along a region given by the interval $[-R,R]$ together with the
half-arc of radius $R$ around the complex origin in the
$\{\Im m>0\}$ half-plane. Since the integrand has degree $-2$ in $t$,
the half-arc contribution vanishes for $R\to\infty$. 

Expand the integrand as a rational function $N(x,y,z,t)/D(x,y,z,t)$
of $x,y,z,t$. The calculation of the discriminant of the (smallest)
denominator polynomial $D(x,y,z,t)=D(t)$, regarded as a
polynomial in $t$, shows that this discriminant has a
non-trivial expansion around $(x,y,z)=(0,0,0)$ (even if
it vanishes in this point). Thus for generic $x,y,z$ of
small norm, $D$ will have only single zeros. These zeros
are explicitly calculable since $D(t)$ decomposes into
quadratic factors in $t$ and $t^2$. Since the solutions
depend continuously on $x,y,z$, to decide which zeros $t_0$ are
relevant, one calculates them for $(x,y,z)=(0,0,0)$. The residues
are then given by $N(t_0)/D'(t_0)$.

The result can be obtained with MATHEMATI\-CA\TM{} \cite{Wolfram}
after some time. It occupied almost 300 lines. Such an expression is
difficult to handle even with the computer. For example, while the
result should have real coefficients, I could not make MATHEMATI\-CA
eliminate the complex units out of it. Nonetheless, substituting small
real values for $x$, $y$ and $z$ showed that $J(x,y,z)$ is indeed real.
After hand manipulation I obtained an expression without occurrences
of $i$, and MATHEMATICA simplified it to about 250 lines. As a check,
expanding $J(x,1,1)$ and $J(x,1,0)$ as power series in $x$ reveals~--
as expected~-- the numbers of all resp.\ $\sg=0$ rational knots. \qed

By applying the same integration to a symmetrized version of $G_2$,
one can also (theoretically) obtain a similar expression for
problem \reference{_typea}.

\begin{rem}
Of course, one could try to solve the integral in \eqref{inl} directly
by residues, without substitution, but it turned out that, when using
MATHEMATICA, manual ``intervention'' was necessary (at least for me)
at an earlier stage. Clearly I tried to avoid this (as long as possible)
with regard to the difficulty of the expressions.
\end{rem}


\begin{rem}
The very useful formula \eqref{intf} seems natural with
harmonic analysis in mind, but 
I could not find, or get referred to, an occurrence of it
in combinatorial literature. M.\ Bousquet-M\'elou pointed out
to me that this product of series is called the Hadamard product.
It has been intensively studied from the point of view of
showing closure properties of certain families of power series
under it (see e.g. \cite{Lipshitz}). A subsequent electronic search
for this term led at least to one reference \cite{Bragg}, where the
integral expression is given explicitly. Thus it appears known in
analysis, even if not popularly. I have no access to that paper and to
the history of the formula, but at least it was discovered independently
by myself. (It occurs also in my previous paper \cite{qeval}.)
\end{rem}

%
%

\section{Further applications\label{Sfa}}

\subsection{\label{Sfg}Applications of the
signature and genus enumeration}

We give another result, in whose proof the integration method
is again applied, and leads to a(n at least electronically) feasible
calculation with a manageably presentable result. (It can be
considered as a special case of proposition \reference{prJ},
up to the different handling of mirror images.)

\begin{corr}
If $c_n$ denotes the number of rational knots of $n$ crossings with
signature $0$ (see last line of table \reference{tab1}), such that
chiral pairs are counted twice, 
then these numbers have a generating function
\[
f_0(x)\,=\,\sum_{n=1}^\infty c_nx^n\,=\,x^4 + 3\,x^6 + 2\,x^7 + 9\,x^8
+ 6\,x^9 + 29\,x^{10} + \dots\,,
\]
which can be expressed in closed form (see figure \reference{tb2}).
Also $\ds\lim\limits_{n\to\infty}\sqrt[n]{c_{n}}=2$.
\end{corr}

\proof The generating function we seek can now be expressed as
\[
\frac{1}{2\pi}\,\li_0^{2\pi}\,G_1\bigl(x,e^{it},e^{-it}\bigr)\,dt\,,
\]
which certainly converges at least for $|x|<\myfrac{1}{2}$. If $G_1$
is a rational function, as in our situation, such an
integral can always be solved. Most generally, with the
standard substitution $z=\tan \myfrac{t}{2}$, it
turns into a rational integral
\[
\,\li_{-\infty}^{\infty}\,G_1\Bigl(x,\frac{1-z^2+2iz}{1+z^2},
\frac{1-z^2-2iz}{1+z^2}\Bigr)\,\frac{dz}{\pi\,(1+z^2)}\,,
\]
which can be solved by the residue method or by partial fraction
decomposition. The result was obtained with MATHEMATI\-CA in a
few minutes. \qed

\begin{rem}
Using the Darboux method (see \cite[\S 5.3]{Wilf}), one can
determine a more precise asymptotic behaviour of the numbers
$c_n$, which is a bit more interesting since their generating
function is not rational. Using the multi-singularity version
of Darboux' theorem \cite[theorem 5.3.2]{Wilf} attributed
to Szeg\"o, and Stirling's formula, one obtains that
the leading term in the asymptotic expansion of $c_n$
is $\ffrac{2^{n-1}}{3\sqrt{2\pi n}}$. (The next
order term contains an oscillating contribution given by a
constant multiple of $(-2)^nn^{-3/2}$.)
\end{rem}

We also remark that we can now easily obtain statistical
data about the distribution of genera and signatures among
rational knots. We give only the (asymptotic) expectation values;
dispersion and the other moments can be determined similarly.

\begin{prop}
The average genus of a rational knot of $n$ crossings behaves
as $n\to\infty$, up to lower order terms, like $\ds\frac{n}{4}$.
The average absolute signature $|\sg|$ behaves like
$\sqrt{\ffrac{2n}{\pi}}$.
\end{prop}

\proof The average genus of a rational knot of $n$ crossings
is given by
\begin{eqn}\label{opr}
\tl g(n)\,:=\,\frac{\ds\sum_{K\in\cC_n}g(K)}{\,|\,\cC_n\,|\,}\,,
\qquad \mbox{with} \qquad
\cC_n\,:=\,\{\,K\ \mbox{rational},\ c(K)=n\ \}\,.
\end{eqn}
Since achiral knots drop exponentially compared to all knots,
it is unimportant for the asymptotics whether we consider knots
up to mirroring or not in $\cC_n$. For convenience, we will assume
for the average genus calculation that we distinguish mirror images,
while for the average absolute signature that we do not.

The behaviour of numerator and denominator in \eqref{opr} are found
by partial fraction decomposition. For the denominator one considers $G_
1(x,1,1)$, and the relevant term one obtains is $\ds\frac{1}{12(1-2x)}$
(which is basically the Ernst-Sumners result). For the numerator
one applies the same procedure to $\ds\frac{\pa G_1}
{\pa y}(x,1,1)$, and finds that $\ds\frac{1}{(1+2x)^2}$ does
not occur and that the coefficient of $\ds\frac{1}{(1-2x)^2}$
is $\myfrac{1}{24}$. Then note that $G_1$ counted in the powers
of $y$ the double genus.

Now consider the average signature (obviously defined).
One calculates $\ffrac{\pa J}{\pa z}(x,1,1)$. The term
whose denominator has zeros of smallest norm is
\[
S(x)\,=\,\frac{(1-x)x^3}{(1+x)(1-2x)^{3/2}\sqrt{1+2x}}\,.
\]
The dominating term thus comes (expectedly) from the zero
$x=\myfrac{1}{2}$. By the Darboux-Szeg\"o theorem, the leading
contribution of this zero is given by
\begin{eqn}\label{ooo}
2^n\es{{n+\myfrac{1}{2}}\choose n}\es\cdot\es
\bigl(S(x)\,\cdot\,(1-2x)^{3/2}\,\bigr)
\raisebox{-0.3em}{$\bigg|_{\,x=\mmyfrac{1}{2}}$}\,.
\end{eqn}
The right factor evaluates to $\ffrac{1}{24\sqrt{2}}$,
and Stirling's formula yields
\[
{n+\myfrac{1}{2}\choose n}\,=\,\frac{\Gm(n+\myfrac{3}{2})}
{\Gm(n+1)\Gm(\myfrac{3}{2})}\,\asymp\,
\frac{\sqrt{n}}{\Gm(\myfrac{3}{2})}\,,
\]
with $\Gm(\myfrac{3}{2})=\ffrac{\sqrt{\pi}}{2}$ and $a_n\,\asymp\,b_n$
meaning $a_n/b_n\to 1$. Then \eqref{ooo} gives 
\[
2^n\,\cdot\,\frac{\sqrt{n}}{12\sqrt{2\pi}}\,,
\]
which, divided by the asymptotical behaviour $2^{n-3}/3$ of the total
number of rational knots up to mirroring, leads to the stated
asymptotics. \qed

\begin{rem}
The generating function $\sum_n\tl g(n)x^n$ of the mean genus $\tl g$
(and also mean $|\sg|$) itself seems difficult to express.
\end{rem}

\begin{rem}
We have for simplicity omitted the following asymptotical terms,
but their contribution is $O(\myfrac{1}{n})$ compared to the
one of the leading term, so that latter alone does not necessarily
give a good approximation. For example, by expanding
$\ffrac{\pa J}{\pa z}(x,1,1)$ as a power series, one finds that
the sum of $|\sg|$ over 1000-crossing rational knots is
about $1.12\times 10^{301}$. Only these first 3 digits are
approximated correctly from the leading term given in the proposition
(when multiplied by the total number of knots).
\end{rem}

\subsection{\label{Sff}Unknotting number one and
the Bleiler conjecture revisited}

Now we return to the enumeration of rational knots of
unknotting number one (with the convention of not distinguishing
mirror images). We promised to give a proof of
corollary \reference{ttt}. For this we consider again
the Conway notation with even numbers, and describe
such notations occurring for unknotting number one knots.

In \cite[\S 3.1]{rtan}, we described exactly arithmetically
which knots $S(p,q)$ give counterexamples to the Bleiler conjecture~--
this occurs iff at least one of the four pairs $(p,\pm q^{\pm 1})$
can be written as $(2mn\pm 1,2n^2)$ with $m>n>1$ coprime, but no one
can be done so such that additionally one of $m$ and $n$ is even.
The main point here is to describe the even-integer notations
for these counterexamples.

\begin{prop}\label{pu1}
Let $K$ be an unknotting number one rational knot, and
$C(a_1,\dots,a_k)$ its Conway notation with non-zero even integers.
Then $(a_1,\dots,a_k)$ is up to transposition of (at least) one of
the following forms:
{\nopagebreak
\def\labelenumi{\arabic{enumi})}\mbox{}\\[-18pt]
\def\theenumi{\arabic{enumi}}
\begin{enumerate}
\item\label{form1} $(a,a_1,\dots,a_l,\pm 2,-a_l,\dots,-a_1)$
with $l\ge 0$ or
\item\label{form2} $(a,a_1,\dots,a_l,\pm 2,a_l',-a_{l-1},\dots,-a_1)$
with $l\ge 1$,
such that $|a_l+a_l'|=2$, and the sign of the absolutely larger one of
$a_l$ and $a_l'$ is opposite to the one of the $\pm 2$ in between.
\end{enumerate}
}
Also, each such sequence is realized by an unknotting number one
rational knot.
\end{prop}

\proof We use the argument in \cite[proof of theorem 1, (ii)
$\Ra$ (iii)]{KanMur}. Take a rational unknotting number one knot
$K=S(2mn\pm 1,2n^2)$ with $(m,n)=1$.
If $n=1$ we have a twist knot, which is of form \reference{form1}.
Thus let $n>1$. Then $m>n$. Kanenobu and Murakami write $m=an+t$,
and now we can choose $a\ne 0$ to be even, possibly having $t<0$.
Then express $n/t$ as a continued fraction. If one of
$n$ and $t$ is even, then one can choose the continued fraction
expression to be only of even integers $a_1,\dots,a_l$, and by the
argument of Kanenobu and Murakami obtains that the form
\reference{form1} can be chosen so that indeed all numbers are even.
Contrarily, every form \reference{form1} clearly gives an unknotting
number one knot.

Now consider the case that both $n$ and $t$ are odd.
Then one can write $n/t$ as a continued fraction, such that all
integers $[[a_1,\dots,a_l]]$ are even except $a_l$, which is odd. One
can also assume that for $l>1$ we have $(a_{l-1},a_l)\ne (\pm 2,\mp 1)$,
and that if $l=1$, then $|a_1|=n>1$. Then use the transformation for
$a>0$
\begin{eqn}\label{yy}
(\dots,a,2,-a,\dots)\,\lra\,(\dots,a+1,-2,1-a,\dots)
\end{eqn}
together with its negated and transposed versions. After an application
of this transformation one obtains a notation of the form
\reference{form2}. 
Also, none of the neighbors of the
middle $\pm 2$ might have become zero under \eqref{yy}, because
we excluded the sequences $(\dots,\pm 2,\mp 1)$ and $(\pm 1)$.
Thus no collapsing occurs, according to the rule \eqref{yyz}.

Finally, it is again easy to see that each sequence of the form
\reference{form2} can be realized.
\qed

Now we can prove corollary \reference{ttt}.

\proof[of corollary \reference{ttt}]
First exclude all twist knots from the consideration. These are
the knots whose notation is of length $2$. They are counted clearly
by
\[
\frac{x^3}{1-x}\,.
\]

Now we count the notations of type \reference{form1} and
\reference{form2} by crossing number. Such notations are unique up
to transposition and negation. To fix the negation ambiguity, we assume
$a>0$.

By similar arguments as before, and using lemma \reference{lerr},
one can find that the generating function of the (remaining,
non-twist-knot) notations of type \reference{form1} by crossing
number is
\[
\frac{2x^6}{(1-x^2-2x^4)(1-x)}\,.
\]
To enumerate type \reference{form2} sequences, just note that
such a sequence of a crossing number $n$ knot bijectively
corresponds to a (non-twist-knot) sequence of a crossing number $n-2$
of type \reference{form1}. Simply raise in latter sequence the absolute
value of one of the neighbors of the middle $\pm 2$ by $2$. The neighbor
is determined by having the opposite sign to the $\pm 2$. Thus
type \reference{form2} sequences are counted by
\[
\frac{2x^8}{(1-x^2-2x^4)(1-x)}\,.
\]

Now we must care about which sequences $w$ are counted several times
up to transposition and possible negation. Clearly $w$ cannot be at
the same time of type \reference{form1} and of type \reference{form2}.
Similarly if both $w$ and $\pm \ol w$ are of type \reference{form1},
or both are of type \reference{form2}, it is easy to see that
$w$ is itself (anti)palindromic, so that it is not generated twice.

Finally, we must care about the case that one of $w$ and $\pm \ol w$
is of type \reference{form1}, and the other one is of type
\reference{form2}. Then one indeed obtains a series of
duplications, namely for the sequences of the form
$(4\, -2)^{k}\,-2\,2\,(-4\, 2)^{k-1}$, and $(2\,-4)^k\,2\,2\,(-2\,4)^k$
with $k\ge 1$. These forms give one knot, in crossing numbers
$8+4r$, $r\ge 0$.

Thus the function we seek is
\[
\frac{x^3}{1-x}\,+\,\frac{2x^6+2x^8}{(1-x^2-2x^4)(1-x)}\,-\,\frac{x^8}
{1-x^4}\,=\,\frac{x^3}{1-x}\,+\,\frac{2x^6}{(1-x)(1-2x^2)}\,-\,
\frac{x^8}{1-x^4}\,,
\]
which is what we claimed. \qed

The proof also gives the following consequence:


\begin{prop}
If $c_n$ is the number
of rational unknotting number one knots of $n$ crossings (chiral pairs
counted once), that do \em{not} provide counterexamples to the
Bleiler conjecture (that is, unknot by one crossing change in their
rational diagrams with all Conway coefficients even), then
\begin{eqnarray*}
\lefteqn{ \sum_{n=1}^{\infty}c_nx^n = \frac{x^3-x^5+2x^6-2x^7}
{(1-x^2-2x^4)(1-x)} } \\ & & \qquad\qquad =\,
x^3+x^4+x^5+3x^6+3x^7+5x^8+5x^9+11x^{10}+11x^{11}+\cdots\,. 
\end{eqnarray*}
\end{prop}

This formula shows that asymptotically $1/3$ of the $n$ crossing
unknotting number one knots do not have the property conjectured by
Bleiler.

\proof This is simply obtained by counting only the twist knots and
the (remaining) ones of type \reference{form1}. \qed


As another consequence we obtain the proof of a weaker form of
Bleiler's conjecture. This form was suggested by, but is
nonetheless still more general than the one proved in \cite{rat2}.

\begin{corr}\label{cBC}
Any unknotting number one counterexample to Bleiler's conjecture has
even leading coefficient $\mc\Dl$ of the Alexander polynomial.
In particular, Bleiler's conjecture holds for unknotting number one
fibered rational knots.
\end{corr}

\proof Use \eqref{pp} and the observation that in type
\reference{form2}, at least one of $a_l$ and $a_l'$ is divisible
by $4$. \qed


The more general version of the fibered Bleiler conjecture
also extends theorem \reference{thFb} to odd values of
$\mc\Dl_K$.
%
%
%
We can, however, obtain a formula even in some cases where the Bleiler
conjecture fails, because we understand well the exceptions.
{}From the proof of theorem \reference{thFb}, and proposition \ref{pu1},
the following can be obtained easily; we leave the proof to the reader.

\begin{prop}
Let $p$ be a square-free positive integer. Then the number of
rational unknotting number one knots $K$ with $\mc\Dl_K=\pm p$ and 
crossing number $n$ is given by
\[
2\left(\,F_{\br{n/2-2-p}}\,+\sum_{r\,:\,r(r+1)\mid p}\,
F_{\br{n/2-1-2r-p/(r+r^2)}}\,\right)\es+\es
\left\{\begin{array}{r@{\es}l}
-1 & \mbox{if}\es (n,p)=(8,2) \\ 1 & \mbox{if}\es n\in\{1+2p,2+2p\} \\
0 & \mbox{otherwise}\end{array}\,\right\}\,.
\]
(In this formula we assume that $r>0$ and that $F_k=0$ if $k<0$.) \qed
\end{prop}

In particular, for square-free odd $p$ and $n\ge 4+2p$ we obtain
$2F_{\br{n/2-2-p}}$, and for $p=2$ and $n\ge 9$ we have
$4F_{\br{n/2-4}}$. When $n\ge 4+2p$, one can use the recursive
behaviour to rewrite the formula also for any other square-free
$p$ to contain only two (mutually index-shifted and bulkily
coefficiented) Fibonacci sequences.
For the remaining, non-square-free values of $p$ one should still obtain
rational generating functions enumerating the corresponding knots,
but these functions will be much less pleasant. (Their shape will
depend on the prime decomposition of the greatest
integer whose square divides $p$.)

Unfortunately, a similarly nice Fibonacci number version is not
possible for achiral unknotting number one knots $K$ of higher
$\big|\mc\Dl_K\big|$, as for each achiral
rational knot $K$ the formula \eqref{pp} shows that $\pm\mc\Dl_K$ is a
square. (In \cite{achir}, we show that this is more generally true
for alternating knots, a result obtained also by Weber and Quach
\cite{VW}.)

\subsection{\label{Sfl}Counting lens spaces}

We conclude our counting results with an application to the
enumeration of lens spaces. In \cite{achir} we gave the
number of different lens spaces of
fundamental group $\bZ_p$. This is equivalent to counting
rational knots by determinant.

\begin{theo}(\cite{achir})
Let $p\ge 3$ be odd.
When considering the lens space $L(p,q)$ and it mirror image
$L(p,-q)=L(p,p-q)$ as equivalent, the number of different
lens spaces with fundamental group $\bZ_p$ is
\begin{eqn}\label{rr}
\frac{1}{4}\left\{\,\phi(p)+r_2^0(p)+2^{\om(p)}\right\}\,,
\end{eqn}
with $r_2^0(p)$ being given by
\[
r_2^0(p)\,=\,\big|\,\{\,(a,b)\in\bN^2\,:\,(a,b)=1,\, a^2+b^2=p\,\}\,\big|\,,
\]
$\om(p)$ denoting the number of different prime divisors of $p$
and $\phi(p)=\,|\bZ_p^*|$ being Euler's totient function.

When distinguishing between $L(p,q)$ and $L(p,-q)$ (if they are
orientation-reversingly inequivalent), the number of
such lens spaces is
\[
\frac{1}{2}\left\{\,\phi(p)+2^{\om(p)}\right\}\,.
\]
\end{theo}

We can now determine the number of lens spaces
which can be obtained by a $p/\pm 2$ surgery along a knot $K$.

\begin{theo}\label{thls}
Let $p\ge 5$ be odd. Then the number $c_p$ of different lens spaces
with fundamental group $\bZ_p$, which are obtainable by a $p/\pm 2$
surgery along a knot $K$, is given by
\begin{eqn}\label{rro}
c_p\,=\,2^{\om((p+1)/2)-1}+2^{\om((p-1)/2)-1}\,+\,
\left\{\begin{array}{r@{\es}l}
-2 & \mbox{if\es $p=p_s$ for some $s\ge 0$}
\\ -1 & \mbox{otherwise}\end{array}\,\right\}\,.
\end{eqn}
In this formula $\om(n)$ denotes as before the number of
different prime divisors of $n$, and
\begin{eqn}\protect\label{rrp}
p_s=\frac{1}{4}\,\left(\,
\bigl( 58 - 41\,{\sqrt{2}} \bigr)\,{{\bigl( 3 - 2\,{\sqrt{2}} \bigr) }^s
} + 
\bigl( 58 + 41\,{\sqrt{2}} \bigr)\,{{\bigl( 3 + 2\,{\sqrt{2}} \bigr) }^s
} \, \right)\,.
\end{eqn}
In this enumeration we consider the lens space $L(p,q)$ and it mirror
image $L(p,p-q)$ as equivalent. If we distinguish them, the number is
\[
2c_p\,-\big|\,\{p\}\,\cap\,\cN\,\big|\,-\,
\big|\,\{p\}\,\cap\,\cS\,\big|\,,
\]
with $\cN\,:=\,\{\,2n^2+2n+1\,:\,n\ge 1\,\}$, $\cS\,:=\,
\{\,q_s\,:\,s\ge 0\,\}$, and
\begin{eqn}\protect\label{rrq}
q_s\,=\frac{1}{3}\,\left( \bigl( 97 - 56\,{\sqrt{3}} \bigr) \,
 {{\bigl( 2 - {\sqrt{3}} \bigr) }^{2\,s}} + 
 \bigl( 97 + 56\,{\sqrt{3}} \bigr)\,
 {{\bigl( 2 + {\sqrt{3}} \bigr) }^{2\,s}}\,+\,1 \right)\,.
\end{eqn}
\end{theo}

\begin{rem}
The numbers $p_s$ and $q_s$ can be given alternatively
in terms of their generating functions 
\[
\sum_{s=0}^\infty p_sx^s\,=\,\frac{29-5x}{1-6x+x^2}\,=\,
29 + 169 x + 985 x^2 + 5741 x^3 + 33461 x^4 + 195025 x^5 + \cdots\,
\]
and
\begin{eqnarray*}
\lefteqn{ \sum_{s=0}^\infty q_sx^s\,=\,
\frac{65-74x+5x^2}{(1-x)(1-14x+x^2)} } \\ & & \qquad\qquad =\,
65 + 901 x + 12545 x^2 + 174725 x^3 + 2433601 x^4 + 33895685 x^5 +
\cdots\,,
\end{eqnarray*}
or by their initial values and linear recursions 
\[
p_s=6p_{s-1}-p_{s-2}\mbox{\es($s\ge 2$)\quad  and\quad}
q_s=15(q_{s-1}-q_{s-2})+q_{s-3}\mbox{\es ($s\ge 3$).}
\]
The $q_s$ do not seem to have been so far of any particular
attention, but the sequence of $p_s$ is listed in \cite{Sloane}
as A001653, with several references.
\end{rem}

\proof[of theorem \reference{thls}]
Let us first prove \eqref{rro}. We know from the arguments of
\cite{KanMur}, which rely on the results of Culler-Gordon-Luecke-Shalen
\cite{CSGL} and Moser \cite{Moser}, that a lens space
$L(p,q)$ is obtainable by $p/\pm 2$
surgery along a knot $K$ if and only if $S(p,q)$ has
unknotting number one. Thus what we claim is equivalent
to counting unknotting number one rational knots (up
to mirroring) by determinant.

It is easy to see, and we remarked it already in \cite{achir},
that, when counting the Schubert notations $S(p,2n^2)$ with $p=
2mn\pm 1$, the first two terms in the formula for $c_p$ just come
from the ways of writing $(p\pm 1)/2=m_{\pm}n_{\pm}$ with
$(m_{\pm},n_{\pm})=1$ up to interchange of $m_{\pm},n_{\pm}$.
We also remarked that the twist knot with determinant $p$ is
counted twice, as occurring in both representations.
The problem was which other duplications occur.

Whenever $n>1$, clearly $n$ determines $m$, and hence $t$. Moreover,
it is easy to see from the proof of proposition \reference{pu1} that
both $n$ and $t$ can be recovered from the forms of both types.
They are just the numerator and denominator of the continued fractions
of $a_1,\dots,a_l$, possibly first undoing the modification
of $a_l$ in type \reference{form2}.
Thus the duplications we sought occur exactly if the even
integer Conway notation, up to reversion and negation,
can be put into these two types of proposition \reference{pu1}
in a different way. But we know now what sequences these are:
we found they belong to one of the two series
$(4\, -2)^{k}\,-2\,2\,(-4\, 2)^{k-1}$
and $(2\,-4)^k\,2\,2\,(-2\,4)^k$ for $k\ge 1$. Also, for any of
these Conway notations exactly two different representations occur.

It is
now a matter of a simple (even if somewhat tedious) calculation
to show that the determinants of corresponding knots are the $p_{s}$.
The easiest way is to note that negated inversion and addition of an
integer correspond to the action of $SL(2,\bZ)$ on the upper half-plane
$\{\Im m>0\}$, given by
\[
\left(\begin{array}{@{\ }c@{\es\,}c@{\ }}
a & b \\ c & d \end{array}\right)\,x\es=\es\frac{ax+b}{cx+d}\,.
\]
Thus the map
\[
[[\,\,\dots,x\,\,]]\ \lmt\ [[\,\,4,-2,\dots,x\,\,]]
\]
can be described by the action of an $SL(2,\bZ)$-matrix.
(Note that prepending a single integer to the iterated fraction,
in our convention, rather than that of \cite{Murasugi2},
cannot be described by such an action because of the sign
switch. However, when prepending two integers, the two sign
changes cancel at the cost of negating the first number prepended.)
This matrix has two distinct Eigenvalues $\lm_{1,2}=3\pm \sqrt{8}$.
Thus for any of the two series the determinants are given by
\[
a\,\lm_1^{2k}\,+\,b\lm_1^k\,+\,c\,+\,d\lm_2^k\,+\,e\lm_2^{2k}\,,
\]
and the coefficients can be determined from the first five values.
Then to verify \eqref{rrp}, one needs to check it only for $s\le 9$.
(Either series are obtained by specifying the parity of $s$.)

When distinguishing $L(p,q)$ and $L(p,-q)$, one needs to
take account of achiral unknotting number one rational knots. We
classified these knots in corollary \reference{crau1} into two series.
(Possibly one can prove the corollary also from the even-integer
notation, but it does not seem worthwhile to get into this now.)
The determinants of the first series are obvious, while those
of the second series $q_s$ are found similarly to $p_s$. \qed

\begin{rem}
One can see that for the doubly counted knots of determinant
$p_s$ in the derivation of \eqref{rro}, one of $2m_+m_-^{\pm 1}$
is a square root of $-1$ in $\bZ_{p_s}^*$. Thus, like the $F_{2k}$,
none of the $p_s$ has a divisor of the form $4r+3$. (This also
follows from the descriptions of the $p_s$ in Sloane's manual.)
\end{rem}

\begin{rem}
In corollary \reference{crau1}, the notation $(1111)$ was
artificially excluded from the second series by writing $(3(12)%
^k1^4(21)^k3)$ instead of $((12)^k1^4(21)^k)$, in order to avoid
mentioning the figure-8-knot twice. Except eventually for its
determinant $5$, it is not clear whether another determinant can
be realized by knots in both series simultaneously, i.e. for some
$s\ge 0$ and $n\ge 1$ we have
\begin{eqn}\label{pee}
q_s\,=\,2n^2+2n+1\,.
\end{eqn}
This problem falls into
the class of polynomial-exponential equations, which have been
for a long time intensively studied and connected to deep
work in number theory (see e.g. \cite{Evertse}). It is
known that, under certain regularity properties (that our example
enjoys), the number of solutions is finite. Apply for instance
theorem 3 of \cite{NP} with $G_m=3q_m-1$ (which form a binary
recurrence with $A=14$ and $B=1$) and $P(x)=6x^2+6x+2$.
While several particular examples have been studied in detail
and some of them solved completely (see \em{loc.\ cit.}\ in
\cite{NP}), ours is apparently not among them, and good
general bounds on the number or size of
solutions are very hard to obtain. 
\end{rem}

At least we have

\begin{prop}
Assume $2n^2+2n+1\in \cS\cap \cN\not=\varnothing$ (with
$\cS$ and $\cN$ defined as in theorem \reference{thls}). Then
{\nopagebreak
\def\labelenumi{\roman{enumi})}\mbox{}\\[-18pt]
\def\theenumi{\roman{enumi})}
\begin{enumerate}
\item\label{yT1} $x=2n+1$ is an integer point on the elliptic curve
\begin{eqn}\label{elc}
y^2\,=\,3x^4+2x^2-5\,=\,(x^2-1)(3x^2+5)\,,
\end{eqn}
with $|x|>3$. ($x=\pm 1,\pm 3$ are obvious points.)
\item\label{yT2} $|\cS\cap \cN|\le 2^{20222}-2\approx
2.68\times 10^{6087}$, and
\begin{eqn}\label{est}
10^{114,000}\,\le\,\min \cS\cap \cN\,\le\,\max \cS\cap \cN\,\le\,
e^{\ds e^{\ds e^{4^{640}5^{16}}}}\,.
\end{eqn}
\end{enumerate}
}
\end{prop}

\proof Consider first part \reference{yT2}. After the most recent
work of Schlickewei and Schmidt \cite{SS,SS3,SS2}, the best 
estimate for the number of solutions $(s,n)$ of \eqref{pee}
one finds is from theorem 2.2(a) of \cite{SS2} applied
on $3q_s-1$ (with $d=t=2$). This gives at most
$2^{20224}$ \em{integer} solutions $(s,n)$. Since we
have the solutions $(-1,1)$, $(-2,0)$, and $n\mt -1-n$ and/or $q\mt
-4-q$ preserve solutions, we obtain at most $2^{20222}-2$
solutions with $s,n\ge 0$.

Using MATHEMATICA, I verified that no solution of \eqref{pee} occurs
for $0\le s\le 10^5$. This establishes the left inequality in
\eqref{est} by evaluating the logarithm of the dominating root
$\log_{10}(2+\sqrt{3})\approx 0.572$. 
To obtain the right inequality, first we transform
the problem to consider integer points on the elliptic curve
\eqref{elc}.

Since the bases $2\pm \sqrt{3}$ appear with even exponents in
\eqref{rrq}, $3q_s-1$ must be an index-$2$-subsequence of a
simpler binary linear recurrence. This recurrence is found to be
\[
\tq_0\,=\,2\,,\quad \tq_1\,=\,4\,,\quad \tq_s\,=\,4\tq_{s-1}-\tq_{s-2}
\,,
\]
and then $3q_s-1\,=\,\tq_{4+2s}$. Define
\[
\tr_0\,=\,0\,,\quad \tr_1\,=\,1\,,\quad \tr_s\,=\,4\tr_{s-1}-\tr_{s-2}
\,.
\]
Then $(2+\sqrt{3})^s\,=\,\ds\frac 12\tq_s+\tr_s\sqrt{3}$. Also
$\ds\frac 14\tq_s^2-3\tr_s^2=1$, because $2+\sqrt{3}$ is a unit of
$\bZ[\sqrt{3}]$ and has norm $1$. Now if
\[
\tq_s\,=\,6n^2+6n+2\,,
\]
then $2\tq_s\,=\,3(2n+1)^2+1\,=\,3x^2+1$, and so
\[
(3x^2+1)^2-48\tr_s^2=16\,,
\]
which yields \eqref{elc} with $y=4\tr_s$. This proves part
\reference{yT1}.

By Baker's work \cite{Baker} the norm $\max(|x|,|y|)$
of an integer solution $(x,y)$ of \eqref{elc} is at most
\[
e^{\ds e^{\ds e^{4^{640}5^{16}}}}\,.
\]
This leads to the upper bound inequality in \eqref{est}, since
$y=4\tr_s\ge \tq_s$ for $s\ge 1$. \qed

\begin{rem}
There have been recently several substantial imrovements of Baker's
result (see e.g. Voutier \cite{Voutier}). However, all these bounds
depend on constants which are (effectively computable but)
not explicitly given. Even if so done, the estimates still seem
too large to close the gap in \eqref{est}. Nonetheless, many special
hyperelliptic equations like \eqref{elc} can be, and have been,
solved completely; this usually requires though, besides use of
general results and computer calculation, a fair bit of
number-theoretically (for me too) advanced extra arguments.
I decided to consult Yu. Bilu about \eqref{elc}; his collaborator
G.\ Hanrot informed me then that, using Magma, he computed that
$x=\pm 1,\pm 3$ are indeed the only solutions.
\end{rem}

\section{Unimodular matrices and the ``fibered'' Bleiler
conjecture\label{SBl}} 

Finally, we make some remarks on how the above discussion on rational
knots can be transferred into a completely arithmetic setting
using iterated fractions. As explained, this is historically
motivated by the result of \cite{rat2}. In particular, we will
see that certain partial cases of the ``fibered'' Bleiler conjecture
follow by purely arithmetic arguments.

Consider $M_{k,l}$ for $k,l\in\bZ$ to be the $2\times2$ matrix 
\[
{\left(\begin{array}{@{}c@{\quad}c@{}}\sgn(k)+4|kl| & 2|k|\sgn(l) \\
2|l| & \sgn(l)\end{array} \right)}\in\Gm_\pm(2)\,.
\]
Here $\Gm_\pm(2)$ denotes the subgroup of $2\times2$ matrices in
$PGL(2,\bZ)$ (that is, matrices of determinant $\pm 1$) with even
lower left entry\footnote{We could avoid the use of `$\sgn$' and
`$|\,.\,|$' and define $M_{k,l}=\left(\begin{array}
{@{}c@{\quad}c@{}}1+4kl & 2k \\ 2l & 1\end{array}\right)$ to be
strictly unimodular, but it is more convenient here to normalize
the matrix so as it to preserve the set of integer vectors $\Bigl\{\,
\left(\begin{array} {@{}c@{}}p \\ q\end{array} \right)\,:\,p>q\ge
1,\,(p,q)=1\,\Bigr\}$.}, and we adopt the convention that $\sgn(0)=1$.
Let $\cM_{Z}$ for $Z\subset \bZ\times\bZ$ be the submonoid
(\em{not} subgroup) of $\Gm_\pm(2)$ generated by $\{\,M_{k,l}\,:\,
(k,l)\in Z\,\}$, and $\cM_Zv=\{\,Mv\,:\,M\in\cM_Z\,\}$ be its
``orbit'' on some vector $v\in\bR^2$.

Then knot theory allows to prove some properties of such orbits.

\begin{prop}\label{pppp}
$\cM_{\bN\times(\bZ\sm\bN_+)}\mybin{1}{0}$ does not contain a vector
of the form $\mybin{2mn\pm 1}{2n^2}$
for any coprime integers $m>n>1$ (and sign choice `$\pm$').
\end{prop}

\proof When identifying $p/q$ for $p>q\ge 1$, $(p,q)=1$
with $\mybin{p}{q}$, then multiplication with $M_{k,l}$ is just the
prepending of $2k,2l$ to the iterated fraction. Then the statement is
just a translation of the fact, proved in \cite{gen1} for arbitrary
knots, that rational positive knots which are not twist knots do not
have unknotting number one. 

Clearly the cases of twist knots are of the given form with $n=1$,
thus we need to pose $n>1$. Now, because of the ambiguity
$S(p,q)=S(p,\pm q^{-1})$ (where the additive and multiplicative
inversions are meant in $\bZ_p^*$), one needs to take
care that for $p=2mn\pm 1$ the even one of the numbers $\sbx{\frac{%
p\pm 1}{2}}$ is not of the form $2n^2$, that is, $4n^2\not\equiv\pm 1
\bmod 2mn\pm 1$. However, for $(m,n)=1$ and $m>n\ge 1$ such a
congruence holds only if $n=1$ (and $m\le 3$), which gets irrelevant
once one poses $n>1$. \qed

A more delicate statement can be made in $\cM_{\bN_+\times(\bZ\sm\{0\})}
$. 

\begin{prop}
Assume $M_{k_1,l_1}\cdot\dots\cdot M_{k_g,l_g}\mybin{1}{0}=
\mybin{2mn\pm 1}{2n^2}$ with $m>n\ge 1$ coprime, $k_i\in\bN_+$ and
$l_i\in\bZ\sm\{0\}$. Then
{\nopagebreak
\def\labelenumi{\roman{enumi})}\mbox{}\\[-18pt]
\def\theenumi{\roman{enumi})}
\begin{enumerate}
\item\label{iT1} $n=g=1$ or $g=2$, and
\item\label{iT2} at most one $l_i$ is negative, and one is exactly if
$2mn\pm 1\equiv 3\bmod 4$ (with the same sign choice as above).
\end{enumerate}
}
\end{prop}

\proof For \reference{iT1}, consider the forms of
proposition \reference{pu1} for the unknotting number one
knot $K=S(2mn\pm 1,2n^2)$. It is easy to see that the
only forms in which all entries of one of the index parities
can be chosen to have all the same sign are form \reference{form1}
for $l\le 1$ and form \reference{form2} for $l=1$; and $g=l+1$.
If $g=1$ (in form \reference{form1}), then $K$ is a twist knot,
so that by the argument in the proof of proposition \reference{pppp},
$n=1$.

For \reference{iT2} compute $\sg(C(2k_1,2l_1,\dots,2k_g,2l_g))$
using lemma \reference{lmsg}, and show that it is given by
$2\#\{\,i\,:\,l_i<0\,\}$. Then use the results of \cite{Murasugi3}
that for any knot $K$, $u(K)\ge |\sg(K)/2|$ and that $|\Dl_K(-1)|-
\sg(K)\equiv 1\bmod 4$. \qed

%
%

While for a number theorist such statements, although possibly not
obvious, may be of insufficient importance, the actual
reason for considering rational knots in this light is because
it may hopefully make the problem of the Bleiler conjecture for fibered
rational knots more arithmetically approachable. In particular,
it can be described by a slight modification of the above propositions.

\begin{prop}\label{cnj1}
$\cM_{\{-1,1\}^{\times2}}\mybin{1}{0}\nowns\mybin{2mn\pm 1}{2n^2}$
for any coprime \em{odd} integers $m>n>1$ and sign choice `$\pm$'.
\end{prop}

\proof
This is a slightly reworded version of Bleiler's conjecture for fibered
rational knots. To explain this, we prove first that
\begin{eqn}\label{M1-1}
\cM_{\{-1,1\}^{\times2}}\mybin{1}{0}\,=\,\left\{
\mybin{p}{q}\,:\,2\mid q,\,p>q>1,\,(p,q)=1,\mbox{ and $S(p,q)$ is
fibered }\right\}\,\cup\,\left\{\mybin{1}{0}\right\}.
\end{eqn}

The fact that $p>q>1$ are coprime and $q$ is even can be verified
by simple arithmetic by virtue of being preserved by multiplication with
any of the $M_{k,l}$. Thus we should explain the fiberedness property.

$S(p,q)$ is fibered for $q$ even iff the iterated fraction of even
integers expressing $p/q$ (which is of even length) contains only
$\pm 2$. When identifying $\mybin{p}{q}$ with $p/q$, then for some
$(\eps_1,\eps_2)\in\{-1,1\}^{\times2}$ the prepending of the two numbers
$2\eps_{1,2}$ to the iterated fraction expression is equivalent
to multiplication with $M_{\eps_1,\eps_2}$ up to change of sign
in one of $p$ and $q$. This discrepancy can be dealt with by
negating all subsequent entries in the indices of the $M$'s
to be multiplied with, which passes the discrepancy through until
it is cancelled with a subsequent sign change, or until the end,
where it gets obsolete. This establishes \eqref{M1-1}.

Now if $(m,n)=1$ with $m>n>1$ odd, then $K=S(2mn\pm 1,2n^2)$
is of unknotting number one, and the proof of proposition
\reference{pu1} shows that its even-integer notation
is of type \reference{form2}. Then we observed that $K$
is not fibered, so that the claim we want to show
follows from \eqref{M1-1}. \qed

Prior to its proof, a computer calculation
verified proposition \reference{cnj1}
for word length $\le 16$, thus in particular Bleiler's conjecture for
rational fibered knots $K$ of crossing number at most 35
(note that the word length in the matrices $M_{k,l}$ is just $g(K)$).
It also showed that vectors of the stated form with both $m$ and $n$
odd, but \em{not}
necessarily coprime occur only for the expected word length 9.
This suggests stronger to conjecture that such vectors can not be
obtained except for word lengths being an odd square.

Note also, that because of the reversal of the
iterated fraction expression with even integers preserves its
form for fibered knots, $\cM_{\{-1,1\}^{\times2}}{1\choose 0}$
contains with ${p\choose q}$ also ${p\choose \pm q^{-1}}$ (where
negation and inversion of $q$ are meant in $\bZ_p^*$
and the sign is chosen so as the number to be even).


One can also try to deduce the above statements in purely arithmetic
terms. As suggested to me by D.~Hejhal, a naive approach is to look at
congruences in the $M_{k,l}$'s. Although this unlikely will lead
to a complete recovering of corollary \reference{cBC} or 
proposition \reference{cnj1}, we record at the end two properties
of genus and determinant that indeed come up this way, namely
from considering $p\bmod 8$ and $q\bmod 4$. 
Working with these congruences, and using
\eqref{pp}, allows, as before, to weaken
the fiberedness assumption on $K$ to $\mc\Dl_K$ being odd, because
each (odd) factor $a_i/2$ of $\mc\Dl_K$ can be made to $\pm 1$ by
adding a multiple of $4$, and this does not affect $a_i\bmod 8$.

\begin{prop}\label{pp2}
If $K$ is an unknotting number one counterexample
to the Bleiler conjecture with $\mc\Dl_K$ odd, then
{\nopagebreak
\def\labelenumi{\roman{enumi})}\mbox{}\\[-18pt]
\def\theenumi{\roman{enumi})}
\begin{enumerate}
\item the genus $g(K)$ of $K$ is also odd,
\item $3\mid g(K)\iff |\Dl_K(-1)| \equiv \pm 1\bmod 8$ and
$3\nmid g(K)\iff |\Dl_K(-1)| \equiv \pm 3\bmod 8$.
\qed
\end{enumerate}
}
\end{prop}

\noindent{\bf Acknowledgement.} I would wish to thank to 
the IH\'ES and its director J.-P. Bourguignon for the invitation
and the hospitality during the period of my visit in
November 1999, when the initial part of this work was done.
At IH\'ES I had the opportunity to speak to several people,
among which the most helpful contribution to the topic
of this paper made D.~Hejhal, L.~Siebenmann and I.~Vardi.
%
I had later several further discussions concerned with knot theory,
combinatorics and number theory related to the subject;
I wish to thank in particular (in approximate chronological
order) T.~Kanenobu, H.~Murakami, A.~Sikora, K.~%
Murasugi, V.~Liskovets, D.~Redmond, J.~Fearnley, Yu.~Bilu, G.~Hanrot
and M.~Bousquet-M\'elou for their helpful remarks and references.
Most of the formulas in sections \S\reference{S6.1} to
\reference{Sfa} have been obtained with the help of MATHEMATICA\TM{}.


{\small
\begin{thebibliography}{CGLS}
\bibitem[Ad]{Adams} C.~C.~Adams, \em{The knot book},
	W.~H.~Freeman \& Co., New York, 1994.
\bibitem[Al]{Alexander} J. W.~Alexander, \em{Topological invariants
	of knots and links}, Trans.\ Amer.\ Math.\ Soc. {\bf 30} (1928),
	275--306.
\bibitem[B]{Baker} A. Baker, \em{Bounds for the solutions of the
	hyperelliptic equation}, 
	Proc. Cambridge Philos. Soc. {\bf 65} (1969), 439--444.
\bibitem[BN]{BarNatanVI} D.~Bar-Natan, {\em On the Vassiliev knot
	invariants}, Topology {\bf 34} (1995), 423--472.
\bibitem[Bl]{Bleiler} S.~A.~Bleiler, \em{A note on unknotting number},
	Math. Proc.~Camb. Phil.~Soc. {\bf 96} (1984), 469--471.
\bibitem[Br]{Bragg} L.\ R.\ Bragg, \em{Trigonometric integrals and
	Hadamard products}, Amer. Math. Monthly {\bf 106(1)}
	(1999), 36--42.
\bibitem[CD]{ChmDuz} S.~V. Chmutov and S.~V. Duzhin, {\em A
	lower bound for the number of Vassiliev knot invariants},
	Topology Appl. {\bf 92(3)} (1999), 201--223. 
\bibitem[Co]{Conway} J.~H.~Conway, \em{On enumeration of knots and
	links}, in ``Computational Problems in abstract algebra''
	(J.~Leech, ed.), 329-358. Pergamon Press, 1969.
\bibitem[Cr]{Cromwell} P. R. Cromwell, {\em Homogeneous links},
	J. London Math. Soc. (series 2) {\bf 39} (1989), 535--552.
\bibitem[CM]{MorCro} \bysame\, and H.~R.~Morton, \em{Positivity of knot
	polynomials on positive links}, J. Knot Theory Ramif. {\bf 1}
	(1992), 203--206.
\bibitem[CGLS]{CSGL} M.\ Culler, C.\ McA.\ Gordon, J.\ Luecke and
	P.\ B.\ Shalen, {\em Dehn surgery on knots},
	Ann. of Math. {\bf 125(2)} (1987), 237--300.
\bibitem[ES]{ErnSum} C. Ernst and D. W. Sumners, {\em
	The Growth of the Number of Prime Knots},
	Proc. Cambridge Phil. Soc. {\bf 102} (1987), 303--315.
\bibitem[Ev]{Evertse} J.-H.\ Evertse, 
	\em{An improvement of the quantitative subspace theorem}, 
	Compositio Math. {\bf 101(3)} (1996), 225--311.
\bibitem[Ga]{Gabai} D.~Gabai, \em{Detecting fibred links in $S^3$},
	Comment. Math. Helv. {\bf 61(4)} (1986), 519--555.
\bibitem[GK]{GK} J.\ R.\ Goldman and L.\ H.\ Kauffman, 
	\em{Rational tangles}, Adv. in Appl. Math.
	{\bf 18(3)} (1997), 300--332. 
\bibitem[Ha]{Haken} W.\ Haken, \em{Theorie der Normalfl\"achen},
	Acta Math. {\bf 105} (1961), 245--375. 
\bibitem[HNK]{HNK} F.\ Hirzebruch, W.\ D.\ Neumann and S.\ S.\ Koh,
	\em{Differentiable manifolds and quadratic forms},
	Lecture Notes in Pure \& Appl.~Math. {\bf 4}, M.~Dekker,
	New York, 1971.
\bibitem[HTW]{HTW} J.~Hoste, M.~Thistlethwaite and J.~Weeks, \em{The
	first 1,701,936 knots}, Math. Intell. {\bf 20 (4)} (1998),
	33--48.
\bibitem[KM]{KanMur} T.~Kanenobu and H.~Murakami, \em{2-bridge
	knots of unknotting number one}, Proc. Amer. Math. Soc. {\bf 
	96(3)} (1986), 499--502.
\bibitem[Ka]{Kauffman} L.~H.~Kauffman, {\em State models and
	the Jones polynomial}, Topology {\bf 26} (1987), 395--407.
\bibitem[KS]{SC} S.\ Kunz-Jacques and G.\ Schaeffer, \em{%
	The asymptotic number of prime alternating links},
	``Formal Power Series and Algebraic Combinatorics'',
	proceedings of the conference in Phenix '01. 
\bibitem[Li]{Lickorish} W. B. R. Lickorish, \em{The unknotting number
	of a classical knot}, in ``Contemporary Mathematics'' {\bf 44}
	(1985), 117--119.
\bibitem[Lp]{Lipshitz} L.~Lipshitz, \em{The diagonal of a $D$-finite
	power series is $D$-finite}, J. Algebra {\bf 113} (1988),
	373--378.
\bibitem[LW]{LW} D.\ Lines and C.\ Weber, \em{N\oe uds rationnels
        fibr\'es alg\'ebriquement cobordants \`a z\'ero}, Topology
	{\bf 22(3)} (1983), 267--283.
\bibitem[MT]{MenThis} W.~W.~Menasco and M.~B.~Thistlethwaite, \em{%
	The Tait flyping conjecture}, Bull. Amer. Math. Soc. {\bf
	25 (2)} (1991), 403--412.
\bibitem[Mo]{Moser} L.\ Moser, \em{Elementary surgery along a torus
	knot}, Pacific J. Math. {\bf 38} (1971), 737--745.
\bibitem[Mu]{Murasugi} K.~Murasugi, \em{Jones polynomial and classical
	conjectures in knot theory}, Topology {\bf 26} (1987), 187--194.
\bibitem[Mu2]{Murasugi2} \bysame, \em{On the braid index of alternating 
	links}, Trans. Amer. Math. Soc. {\bf 326 (1)} (1991), 237--260.
\bibitem[Mu3]{Murasugi3} \bysame\,, \em{On a certain numerical
	invariant of link types}, Trans. Amer. Math. Soc. {\bf 117}
	(1965), 387--422.
\bibitem[Mu4]{Murasugi4} \bysame\,, \em{On a certain subgroup of
	the group of an alternating link}, Amer. J. Math. {\bf 85}
	(1963), 544--550.
\bibitem[N]{Nakamura} T.~Nakamura, \em{Positive alternating links
	are positively alternating}, J. Knot Theory Ramifications
	{\bf 9(1)} (2000), 107--112.
\bibitem[NP]{NP} I. Nemes and A. Peth\"o, \em{Polynomial values in
	linear recurrences II},
	J. Number Theory {\bf 24(1)} (1986), 47--53. 
\bibitem[Ro]{Rolfsen} D.~Rolfsen, {\em Knots and links} (second
	edition), Publish or Perish, 1993.
\bibitem[SS]{SS} H. P. Schlickewei and W. M. Schmidt,
	\em{On polynomial-exponential equations},
	Math. Ann. {\bf 296(2)} (1993), 339--361.
\bibitem[SS2]{SS3} \bysame\,and \bysame, \em{The intersection of
	recurrence sequences}, Acta Arith. {\bf 72(1)} (1995), 1--44.
\bibitem[SS3]{SS2} \bysame\,and \bysame, \em{The number of solutions
	of polynomial-exponentials equations},
	Compositio Math. {\bf 120(2)} (2000), 193--225.
\bibitem[Sh]{Schubert} H.\ Schubert, \em{Knoten mit zwei Br\"ucken},
	Math.\ Z. {\bf 65} (1956), 133--170.
\bibitem[Sl]{Sloane} N.~J.~A.~Sloane, {\em The On-Line Encyclopedia of
	Integer Sequences}, accessible at the web address
	{\verb|http://www.research.att.com/~njas/sequences/index.html|}.
\bibitem[St]{positive} A.~Stoimenow, {\em Positive knots, closed braids,
	and the Jones polynomial}, preprint {\tt math/9805078}.  
\bibitem[St2]{gen1} \bysame, {\em Knots of genus one},
	Proc. Amer. Math. Soc. {\bf 129(7)} (2001), 2141--2156.
\bibitem[St3]{restr} \bysame, {\em On some restrictions to the values 
	of the Jones polynomial}, preprint.
\bibitem[St4]{rtan} \bysame, \em{On unknotting numbers and knot
	trivadjacency}, partly joint work with N.~Askitas, preprint.
\bibitem[St5]{rat2} \bysame, \em{Fibonacci numbers and 
	the ``fibered'' Bleiler conjecture}, Int.\ Math.\ Res. Notices
	{\bf 23} (2000), 1207--1212.
\bibitem[St6]{nlpol} \bysame, {\em On the number of links and link
	polynomials}, preprint.
\bibitem[St7]{achir} \bysame, {\em Square numbers, spanning trees and
	invariants of achiral knots}, preprint {\tt math.GT/0003172}.
\bibitem[St8]{qeval} \bysame, {\em Branched cover homology and $Q$
	evaluations}, Osaka J. Math. {\bf 39(1)} (2002), 13--21.
\bibitem[SV]{SV} \bysame\ and A.\ Vdovina, {\em Counting
	alternating knots by genus}, preprint.
\bibitem[ST]{SunThis} C.~Sundberg and M.~B.~Thistlethwaite, {\em The
	rate of growth of the number of prime alternating links and
	tangles}, Pacific Journal of Math. {\bf 182 (2)} (1998),
	329--358.
\bibitem[Th]{This} M.~B.~Thistlethwaite, {\em A spanning tree
	expansion for the Jones polynomial}, Topology {\bf 26} (1987),
	297--309.
\bibitem[VW]{VW} Quach H. C\^am V\^an and C. Weber, {\em On the
	topological invariance of Murasugi special components of an
	alternating link}, preprint.
\bibitem[V]{Voutier} Paul M. Voutier, \em{An upper bound for the size
	of integral solutions to $Y\sp m=f(X)$},
	J. Number Theory {\bf 53(2)} (1995), 247--271.
\bibitem[W]{Welsh} D. J. A. Welsh, \em{On the number of knots and
	links}, Sets, graphs and numbers (Budapest, 1991), 713--718,
	Colloq. Math. Soc. J\'anos Bolyai, {\bf 60}, North-Holland,
	Amsterdam, 1992.
\bibitem[Wi]{Wilf} H. S. Wilf, \em{Generatingfunctionology},
	Academic Press, Inc., Boston, MA, 1990.
\bibitem[Wo]{Wolfram} S. Wolfram, {\em Mathematica --- a system
   	for doing mathematics by computer}, Addison-Wesley, 1989.
\bibitem[Yo]{Yokota} Y.\ Yokota, {\em Polynomial invariants of
	positive links}, Topology {\bf 31(4)} (1992), 805--811.
\bibitem[Za]{Zagier} D. Zagier, \em{Vassiliev invariants and a
	strange identity related to the Dedekind eta-function},
	Topology {\bf 40(5)} (2001), 945--960.
\bibitem[ZZ]{ZZ} P.\ Zinn-Justin and J.-B. Zuber, \em{On the counting
	of colored tangles}, J. Knot Theory Ramifications
   	{\bf 9(8)} (2000), 1127--1141.
\end{thebibliography}
}
\end{document}